\documentclass[11pt]{article} 

\usepackage[utf8]{inputenc} 

\usepackage{geometry} 
\usepackage{graphicx} 

\usepackage{booktabs} 
\usepackage{array} 
\usepackage{paralist} 
\usepackage{verbatim} 
\usepackage{subfig} 
\usepackage{amsmath,amsthm,amssymb,color,epic,eepic,cite}
\usepackage{hyperref}
\usepackage{tikz}
\usetikzlibrary{matrix}
\usepackage{fancyhdr} 
\usepackage{sectsty}
\allsectionsfont{\sffamily\mdseries\upshape} 

\usepackage[nottoc,notlof,notlot]{tocbibind} 
\usepackage[titles,subfigure]{tocloft} 

\newcommand{\Z}{{\Bbb Z}} 
\newcommand{\C}{{\Bbb C}} 
\newcommand{\K}{{\Bbb K}}

\newcommand{\hbs}{{\bar{H}^*}}
\newcommand{\tpt}{{\tilde{\otimes}}} 
\newcommand{\D}{{\mathcal{D}}}

\newcommand{\cD}{{\cal D}}

\newcommand{\cR}{{\cal R}}

\newcommand{\hL}{\widehat{L}}

\newcommand{\la}{\lambda}

\newcommand{\al}{\alpha}
\newcommand{\e}{\epsilon}
\newcommand{\vep}{\varepsilon}

\newcommand{\bH}{\bar{H}}

\newcommand{\hOmega}{\widehat{\Omega}}

\newcommand{\hf}{\widehat{f}}

\newcommand{\hV}{\widehat{V}}
\newcommand{\hW}{\widehat{W}}

\newcommand{\noi}{{\noindent}}
\newcommand{\nn}{{\nonumber}}
\newcommand{\bea}{\begin{eqnarray}}
\newcommand{\ena}{\end{eqnarray}}
\newcommand{\beit}{\begin{itemize}}
\newcommand{\enit}{\end{itemize}}

\newcommand{\be}{\begin{eqnarray*}}
\newcommand{\en}{\end{eqnarray*}}
\newcommand{\lb}[1]{\label{#1}}
\newcommand{\End}{{\rm End}}

\def\infq4p#1{{(#1;q^4,p)_\infty}}

\newcommand{\tot}{\widetilde{\otimes}}

\newcommand{\mmatrix}[1]{\begin{matrix} #1 \end{matrix}}

\font\teneufm=eufm10
\font\seveneufm=eufm7
\font\fiveeufm=eufm5
\newfam\eufmfam
\textfont\eufmfam=\teneufm
\scriptfont\eufmfam=\seveneufm
\scriptscriptfont\eufmfam=\fiveeufm
\def\frak#1{{\fam\eufmfam\relax#1}}
\let\goth\frak
\newcommand{\slth}{\widehat{\goth{sl}}_2}
\newcommand{\glth}{\widehat{\goth{gl}}_2}
\newcommand{\slt}{\goth{sl}_2}

\newcommand{\U}{{\mathcal{U}}}
\newcommand{\LL}{{\cal L}} 

\newcommand{\hh}{\goth{h}}
\newcommand{\h}{H}

\font\seventeeneufm=eufm10 scaled\magstep3   
   
\newcommand{\slthBig}{\widehat{\mbox{\seventeeneufm sl}}_2} 
\newcommand{\fp}{{F^+(u_1)}}

\newcommand{\Ep}{{E^+(u_1)}}

\newcommand{\fpt}{{F^+(u_2)}}

\newcommand{\Ept}{{E^+(u_2)}}

\newcommand{\kpone}{{K_1^+(u_1)}}

\makeatletter
\@addtoreset{equation}{section}
\makeatother

\newtheorem{thm}{Theorem}[section]
\newtheorem{prop}[thm]{Proposition}

\newtheorem{cor}[thm]{Corollary}

\newtheorem{dfn}[thm]{Definition}

\title{Dynamical FRT construction of $U_{q,x}(gl_N)$}
\author{Bharath Narayanan\\
The Pennsylvania State University\\
narayana@math.psu.edu}

\begin{document}
\bibliographystyle{unsrt}

\vspace{2cm}
\begin{center}
{\Large \bf Representations of the Dynamical Affine Quantum Group $U_{q,x}(\slth)=U_{q,\lambda}(\slth)$  and Hypergeometric Functions.\\[10mm]}
 
{\large  Bharath Narayanan}\\[6mm]
{\it Department of Mathematics, 
\\The Pennsylvania State University\\
       narayana@math.psu.edu}\\[10mm]
\end{center}

\begin{abstract}
\noindent 
The representation theory of the Hopf algebroid $U_{q,x}(\slth)=U_{q,\la}(\slth)$ is initiated and it is established that the intertwiner between the tensor products of dynamical evaluation modules is a well-poised balanced $_{10}W_9$ symbol, confirming a conjecture of H.Konno, that the degeneration of the elliptic $_{12}V_{11}$ series to  $_{10}W_9$ can be proven based on the representation theory of  $U_{q,x}(\slth)$, viewed as the degeneration of the elliptic algebra $U_{q,p}(\slth)$ as $p \to 0$.

\end{abstract}
 
\section{Introduction}
\begin{figure}[here]
\includegraphics[width=1.01\textwidth]{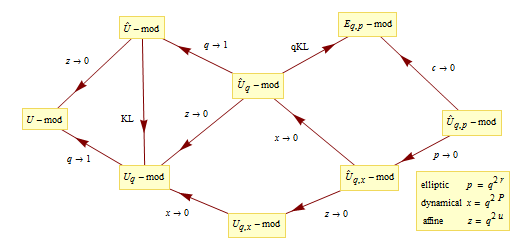}
\caption{Finite-Dimensional Modules, Degenerations and qKazhdan-Lusztig Functors.}
\label{fig:figure1.1}
\end{figure}

The representation theory of quantum algebras is one of the pillars of modern mathematics.   The interplay between the finite-dimensional irreducible modules over the quantum, quantum affine and dynamical quantum affine algebras is illustrated in Figure \ref{fig:figure1.1}.  One of the major achievements is the (quantum) Kazhdan-Lusztig functor
and another is the elliptic quantum group of H.Konno.  The advantage with the approach used in the latter is that one can easily derive both the finite and the infinte dimensional representations from those for $U_q(\slth)$, which are well known, and the coalgebra structure is easily accesible,  facilitating the study of  tensor products of dynamical representations of the associated $H$-Hopf algebroids.  
The Kazhdan-Lusztig functor and quantum-Kazhdan-Lusztig functor \cite{EM} are isomorphisms of tensor categories of finite-dimensional representations:
\be
 U(\slth)-mod\ \simeq\ U_q(\slt)-mod  \quad \text{ and } \quad U_q(\slth)-mod\ \simeq \ E_{q,p}(\slt)-mod, \en 
respectively, as shown in \ref{fig:figure1.1}.

\noi


(i)
$E_{q,p}=E_{\tau,\eta}(\slt)$ is the elliptic quantum group  of Felder \cite{Felder}.  Its representation theory is initiated in \cite{FV} ($R$-matrix representations) and \cite{KNR} (modules and comodules).  In \cite{KNR} it is shown that the intertwiners of comodules are given by a  ${}_{10}W_{9}$ symbol.   

(ii)
$\widehat{U}_{q,p}=U_{q,p}(\slth)$ is the elliptic quantum group of Jimbo, Konno, et al \cite{JKOS2, K2}.  The representation theory is examined in \cite{K2} and a ${}_{12}V_{11}$ symbol is derived using intertwiners of evaluation modules (these modules coincide with (i) when the central element $c=0$.)

(iii)
$U_{q,x}=F_{q,\la}{[{SL2}]}$ is the trignometric dynamical quantum group of Etingof-Varchenko \cite{EV1,EV2} and P. Xu \cite{ping}.  Its dynamical modules, comodules and $R$-matrix representations are investigated in \cite{Rosengren}.  On a set-theoretic level, an equivalence of representation categories of the FRT-bialgebroid $A(R)$ and the trignometric $R$-matrix $R(\la)$ is proven in \cite{ST}.

(iv)
$\widehat{U}_q=U_q(\slth)$ is the usual quantum affine algebra.  Its finite-dimensional modules are essentially the same as modules over $\widehat{U}=U(\slth)$.  We review the construction of $\widehat{U}_q$-modules \cite{CP} in Section 3.1.

(v) 
$U_q=U_q(\slt)$ is the quantum group of Drinfeld-Jimbo.  It is well-know that its finite-dimensional representations are similar to representations of $U=U(\slt)$ \cite{CPbook}.  Essentially, $U$-mod and $U_q$-mod are isomorphic as categories but not as tensor categories.

(vi) $U_{q,x}=U_{q,x}(\slth)$ is constructed in \cite{me1} and its representation theory is the topic of this treatise. A definition based on $RLL$-relations is suggested, but not expounded, in \cite{clad}.
 The non-elliptic, non-dynamical and non-affine degenerations can be considered on the level of $R$-matrices,  algebras,  modules and intertwiners ($6j$-symbols).  
In this article, we consider the dynamical algebra $U_{q,x}(\slth)$, defined in terms of Drinfeld currents and a Heisenberg algebra $\left\{P,Q\right\}$, or equivalently by an FRT construction using $RLL$ relations \cite{me1}.  The subalgebras $U^+_{q,x}(\slth)$ (resp  $U^-_{q,x}(\slth)$) are defined using the $L^+$-operator (resp. $L^-$-operator) only and coincide with the subalegbras of non-positive (resp. non-negative) powers  of $z$ for the elliptic quantum group  $U_{q,p}(\slth)$ of H. Konno when $p\rightarrow0$. 

 We study the dynamical representations of $U_{q,x}(\slth)$, specifically the pseudo highest-weight modules, and prove a criterion for finite-dimensionality using dynamical Drinfeld polynomials.  The main examples are the evaluation modules.  Using the intertwiners of the tensor products of these modules, we confirm a conjecture of H. Konno \cite{K2} that the $_{10}W_9$ series can be obtained based on the representation theory of $U_{q,x}(\slth)$.  To be precise, we show that the elliptic 6j symbol $_{12}V_{11}$ degenerates into Wilson's $_{10}W_9$ hypergeometric series (trigonometric 6j symbol) as $p \to 0$, analogous to the corresponding algebra degeneration $U_{q,p}(\slth)$ to $U_{q,x}(\slth)$.  The representations constructed here are consistent with the $R$-matrix representations.
We explicitly derive the expressions for the action of the positive half-currents, negative half-currents and total currents as well as of the $L^\pm(u)$ on the dynamical evaluation modules and confirm that they are consistent with their non-dynamical degenerations on the one hand, while they are they coincide with the non-elliptic degenerations of the corresponding $U_{q,p}$-modules, on the other.

The fusion or intertwining operators are of primary interest from a physical viewpoint, in the models associated to spin recoupling, as well as 
 for mathematical purposes - hypergeometric functions, category theory and non-commutative geometry.  The dynamical quantum Yang-Baxter appears first in the work of Wigner (1940) \cite{W}.  The $6j$ symbols have a deep history \cite{Rosengren6j} and the generalization to infinite-dimensional dynamical quantum algebras is an area of active research \cite{K2}.  

{\em Outline of this article.} Section 2 reviews the constructions of $U_{q,x}(\slth)$ and $U(R)$ using dynamical Drinfeld currents and using the $L$-operators, respectively as given in \cite{me1}.  The H-Hopf algebroid structure is discussed in this section.  Section 3 begins with a review of the representation theory of $U_q(\slth)$ and of dynamical representations of $H$-algebras.   Next, we define the pseudo-highest weight modules of $U_{q,x}(\slth)$ and prove a criterion for their finiteness using dynamical Drinfeld polynomials.  The main example of evaluation modules in terms of the different realizations is presented next.  The tensor products of dynamical representations are defined in the final section where we prove a criterion for them to be pseudo highest-weight.  The proof of Konno's conjecture appears at the end of the section.  The (positive and negative) dynamical $RLL$ relations are expanded in Appendix A.

\section{The Construction of the $H$-algebra \texorpdfstring{$\U:=U_{q,x}(\slth)$}{Uqx}}
\subsection {Standard Drinfeld Realization of \texorpdfstring{$U_{q}(\slthBig)$}{quantum affine sl2}}
Consider a field $\K \supseteq \C$. The following standard presentation, adapted from \cite{JM} by replacing $c$ with $-c$, first appeared in  \cite{Drinfeld}.

\begin{dfn}[Drinfeld]\lb{defstandard}
$\K[U_q(\slth)]$ is the associative algebra over $\K$ 
generated by 
$a_n\ (n\in \Z_{\not=0})$, $x_n^\pm\ (n\in \Z)$, 
$h$,  $c$ and $d$, with the defining relations 
\be
&&c :\hbox{ central },\nn\\
&& [h,d]=0,\quad [d,a_{n}]=n a_{n},\quad 
[d,x^{\pm}_{n}]=n x^{\pm}_{n}, \nn\\
&&[h,a_{n}]=0,\qquad [h, x_n^\pm]=\pm 2 x_n^\pm,\nn\\
&&
[a_{m},a_{n}]=\frac{[2n]_{q}[c n]_{q}}{n}\delta_{n+m,0},\nn\\
&&[a_{n},x_m^+]=\frac{[2n]_{q}}{n}q^{\frac{c|n|}{2}}x_{m+n}^+,\nn\\
&&[a_{n},x_m^-]=-\frac{[2n]_{q}}{n}q^{-\frac{c|n|}{2}} x_{m+n}^-,\nn\\
&&x_{m+1}^\pm x_n^\pm - q^{\pm 2}x_n^\pm x_{m+1}^\pm =  q^{\pm 2}x_m^\pm x_{n+1}^\pm - x_{n+1}^\pm x_m^\pm ,\nn\\
&&[x^+_n,x^-_m] =\frac{1}{q-q^{-1}}\left(q^{\frac{c(m-n)}{2}}\psi _{m+n}-q^{\frac{c(n-m)}{2}}\varphi _{m+n}\right). \nn
\en

Denote $x^\pm(z)=\sum_{n\in \Z}x^\pm_{n} z^{-n}$, and
the auxillary currents $\psi_n$,\ $\varphi_{-n}$, ($n \geq 0$) by
\be
&&\sum_{n\geq 0}\psi_n z^{-n}=q^{h}
\exp\left( (q-q^{-1}) \sum_{n>0} a_{n}z^{- n}\right),\quad 
\sum_{n\geq 0}\varphi_{-n}z^n=q^{-h}
\exp\left(-(q-q^{-1})\sum_{n>0} a_{-n}z^{ n}\right).\\
\en

\end{dfn}

\subsection{Definition of the $H$-algebra \texorpdfstring{$U_{q,x}(\slth)$}{Uqx}}\lb{subsectionmaindef}

The $H$-algebra structure on  $U_{q,x}(\slth)$ is defined in \cite{me1}.  Let us review the construction.
Define a Heisenberg algebra $\mathbb{H}$ with generators $P$ and $Q$ such that $[P,Q]=-1$. Denote $\mathbb{H}=\text{$\mathbb{C}$Q}\oplus \text{$\mathbb{C}$P}= H^{* }\oplus  H$ with the pairing given by $<Q,P> = 1$ and $<x,y> =0$ for all other $x,y$ (for example, we can choose $Q=\frac{\partial }{\partial P}$).
 Let $\hbs = \text{$\mathbb{Z}$Q}$. 
  Consider the isomorphism  $\Phi :\mathcal{Q} \to  \bar{H}^*,e^{\text{n$\alpha $}_1} \mapsto \e^{\text{nQ}}$.  We will identify $\hbs$ with its group algebra $\mathbb{C}[\hbs]$ by $\alpha  \mapsto \e^{\alpha }$.  Just like in the elliptic case considered in \cite{K2}, we identify
$\widehat{f}=f(P) \in\C [H]$ and  meromorphic functions on $\h^*$ by 
\be
\widehat{f}(\mu)=f<\mu,P>,\quad \mu\in {H}^*
\en
and recover the field of  meromorphic functions  $M_{{\h}^*}$ on ${\h}^*$ 
\be
&&{M}_{{\h}^*}=\left\{ \widehat{f}:{\h}^*\to \C\ \left|\ 
\widehat{f}=f(P)\in \C [H]\right.\right\}.
\en
Let $\K:=\C[H]$ and define the $H$-algebra $U_{q,x}(\slth):=\K[U_q{(\slth)}] \otimes \C[\bar{H}^*]$.
The moment maps are given by:
\bea
&&\mu_l(\hf)=f(P+h), \qquad   \mu_r(\hf)=f(P). 
\ena
The $H$-bigrading is:
\be
&&U_{q,x}(\slth)=\bigoplus_{\alpha,\beta \in H^*}{U_{q,x}(\slth)_{\alpha,\beta}},\\
&& U_{q,x}(\slth)_{\alpha,\beta}=
\left\{\ x\in U_{q,x}(\slth) \left|\ \mmatrix{q^{P+h}xq^{-(P+h)}=q^{<\alpha, P>}x, \cr
q^{P}xq^{-P}=q^{<\beta, P>}x\cr}\ \right.\right\}.
\en

\noi
For $a,b$ in $\mathbb{C}[U_q(\glth)]$, the multiplication in  $U_{q,x}(\slth)$ is defined through the expression
\bea \lb{mult}
&&\left(f(P)a \otimes  e^{\alpha }\right)\cdot \left(g(P)b \otimes  e^{\beta }\right)=f(P)g(P+<\alpha ,P>)\text{ab}\otimes e^{\alpha +\beta }.
\ena

\noi
Let $\K=\C[H]$ and $\widehat U_{q,x}(\slth):=\K[U_q(\slth)] \otimes \C[\hbs]$ and define the operators {\footnote{These were obtained by suitable degenerations of the expressions for the elliptic $K(z)$ and $H^\pm(u)$ given in eqns (3.15, 3.25, 3.29) in \cite{JKOS2} (without the elliptic shift by the central element: $p \to pq^{-2c}$).}}
\vspace{-4pt}
\bea \lb{defnKplus}
&& K^+(u)=\exp\left(\sum_{n > 0} \frac{[n]_q}{[2n]_q}(q-q^{-1})a_n q^{-(2u+1)n}\right)  e^Q q^{\frac{h}{2}},\\
&& K^-(u)=\exp\left(-\sum_{n > 0} \frac{[n]_q}{[2n]_q}(q-q^{-1})a_{-n} q^{(2u+1)n}\right)  e^Q q^{-\frac{h}{2}}.
\ena
\begin{dfn}[Dynamical Currents] We define 
\bea
&&  {K_1^\pm}(u)=K^{\pm}(u-1), \quad  {K_2^\pm}(u)=K^{\pm}(u)^{-1}, \\
&& {E}(u)=x^+(z)e^{2Q}, \quad  {F}(u)=x^-(z), \quad  {H^{\pm }}(u)= {K_1^{\pm }}(u) {K_2^{\pm }}(u)^{-1}. \lb{hatefhdefn}
\ena
\end{dfn}
Define the derived subalgebra  $U_{q,x}(\slth)'$ of $U_{q,x}(\slth)$ by the same relations as Definition \ref{defstandard}, but without the element $d$.  The next proposition describes the commutation relations for  $U_{q,x}(\slth)'$.  We note that they coincide with the elliptic relations in \cite {JKOS2} at $p=0$.  The functions $\rho^\pm(u)$ are given in \eqref{rho} and \eqref{rhominus}.
\begin{prop}\lb{Defrelns}. Let $\eta (x) = 1-q^{2x}.$  The generators of $U_{q,x}(\slth)$ satisfy:
\be
&&q^c:\hbox{ {\rm central}}, \lb{u1}\\
&& [h,E(u)]=2E(u),\quad [h,F(u)]={-2}F(u),\\
&&K_1^{\pm }(u)K_1^{\pm }(v)=K_1^{\pm }(v)K_1^{\pm }(u), \lb{u3} \quad
K_2^{\pm }(u)K_2^{\pm }(v)=K_2^{\pm }(v)K_2^{\pm }(u), \lb{u4}
\quad
K_1^{\pm }(u)K_2^{\pm }(v)=K_2^{\pm }(v)K_1^{\pm }(u),
\lb{u2}\\
&&K_1^+(u)K_1^-(v)=\frac{\rho^{+} \left(u-v-\frac{c}{2}\right)}{\rho^{+} \left(u-v+\frac{c}{2}\right)}K_1^-(v)K_1^+(u),\lb{u5}\quad
K_2^+(u)K_2^-(v)=\frac{\rho^{+} \left(u-v-\frac{c}{2}\right)}{\rho^{+} \left(u-v+\frac{c}{2}\right)}K_2^-(v)K_2^+(u),\lb{u7}\\
&&K_1^+(u)K_2^-(v)=\frac{\rho^{+} \left(u-v-\frac{c}{2}\right)}{\rho^{+} \left(u-v+\frac{c}{2}\right)}\frac{\eta \left(u-v+\frac{c}{2}\right)\eta \left(u-v-\frac{c}{2}-1\right)}{\eta \left(u-v-\frac{c}{2}\right)\eta \left(u-v+\frac{c}{2}-1\right)}K_2^-(v)K_1^+(u),\lb{u9}\\
&&K_2^+(u)K_1^-(v)=\frac{\rho^{+} \left(u-v-\frac{c}{2}\right)}{\rho^{+} \left(u-v+\frac{c}{2}\right)}\frac{\eta \left(u-v-\frac{c}{2}\right)\eta \left(u-v+\frac{c}{2}+1\right)}{\eta \left(u-v+\frac{c}{2}\right)\eta \left(u-v-\frac{c}{2}+1\right)}K_1^-(v)K_2^+(u),\lb{u99}\\
&&K_1^\pm(u)^{-1}E(v)K_1^\pm(u)=\frac{q\eta \left(u-v\mp\frac{c}{4}-1\right)}{\eta \left(u-v\mp\frac{c}{4}\right)}E(v),\lb{u10}\quad\
K_1^\pm(u)F(v)K_1^\pm(u)^{-1}=\frac{q\eta \left(u-v\pm\frac{c}{4}-1\right)}{\eta \left(u-v\pm\frac{c}{4}\right)}F(v),\lb{u12}\\
&&K_2^\pm(u)^{-1}E(v)K_2^\pm(u)=\frac{q^{-1}\eta \left(u-v\mp\frac{c}{4}+1\right)}{\eta \left(u-v\mp\frac{c}{4}\right)}E(v),\lb{u14}\
K_2^\pm(u)F(v)K_2^\pm(u)^{-1}=\frac{q^{-1}\eta \left(u-v\pm\frac{c}{4}+1\right)}{\eta \left(u-v\pm\frac{c}{4}\right)}F(v),\lb{u16}\\
&&E(u)E(v)=\frac{q^{-2}\eta\left(u-v+1\right)}{\eta \left(u-v-1\right)}E(v)E(u) \lb{u18}, \quad
F(u)F(v)=\frac{q^2\eta \left(u-v-1\right)}{\eta\left(u-v+1\right)}F(v)F(u),\lb{u19}\\
&&\left[E(u),F(v)\right]= \frac{1}{q-q^{-1}}\left(\delta \left(q^{-c}\frac{z}{w}\right)H ^+\left(q^{\frac{c}{2}}w\right)-\delta \left(q^c\frac{z}{w}\right) H^-\left(q^{\frac{-c}{2}}w\right)\right), \text{ where }  \delta(t)=\sum_{n \in \Z}t^n.\lb{u20}\\
\en
\end{prop}
\subsection{Half Currents}
\begin{dfn}[Positive Half-Currents]\lb{halfcurrentintegrals}
\be
&&E^+(u)
=q^{-1}a_1 \oint_{C_1} E(u') 
\frac{q^{2(P-1)}\eta\left(u-u'-c/4-P+1\right)\eta(1)}
{\eta(u-u'-c/4)\eta(P-1)}
\frac{dz'}{2\pi i z'}\ ,
\lb{Eplus}\\
\en
\be
&&F^+(u)
=q^{-1}a_2 \oint_{C_2} F(u') 
\frac{\eta(u-u'+c/4+P+h-1)\eta(1)}{\eta(u-u'+c/4)\eta(P+h-1)}
\frac{dz'}{2\pi i z'}. 
\lb{Fplus}
\en

\noi
Here the contours are chosen as
\be
&&C_1 : 0<|z'|<|q^{c/2}z|, \qquad
C_2 : 0<|z'|<|q^{-c/2}z|,
\lb{C}
\en
and the constants $a_1,a_2$ are chosen to satisfy 
$q^{-2}a_1a_2\eta (1)^2 = -1.$

\end{dfn}
\noi
We will use the following Laurent expansion valid in the domain $|z|>1$:
\bea \lb{etaus}
&& \frac{\eta(u+s)}{\eta(u)\eta(s)}=-\sum_{n\in \Z_{\ge 0}}\frac{1}{1-q^{-2s}\delta_{0,n}}z^{-n}.
\ena
This implies that 
\bea \lb{etadelta}
&&\frac{\eta(u+s)}{\eta(u)\eta(s)}+\frac{\eta(-u-s)}{\eta(-u)\eta(-s)}=-\delta(q^{2u}),
\ena
with the first summand expanded around $z=\infty$ and the second, around $z=0$.
\noi
Definition  \ref{halfcurrentintegrals} of the positive half-currents is restated and extended to the negative half-currents, using \eqref{etaus}.
\begin{dfn}[The Series Definition of Half-Currents] \lb{halfcurrentseries}
\be
&&E^+(u)=
e^{2Q}q^{-1}a_1\eta(1)\left(e_0\frac{1}{1-q^{-2(1-P)}}+ \sum _{n>0} e_n(q^{-\frac{c}{2}}z)^{-n} \right),
\lb{Epluss} \\
&&E^-(u)
=-e^{2Q}q^{-1}a_1\eta(1)\left(e_0\frac{1}{1-q^{2(1-P)}}+ \sum _{n<0} e_n(q^{\frac{c}{2}}z)^{-n}\right),
\lb{Eminuss}\\
&&F^+(u)=
-q^{-1}a_2\eta(1)\left(f_0\frac{1}{1-q^{-2(P+h-1)}}+ \sum _{n>0} f_n(q^{\frac{c}{2}}z)^{-n} \right),
\lb{Fpluss} \\
&&F^-(u)
=q^{-1}a_2\eta(1)\left(f_0\frac{1}{1-q^{2(P+h-1)}}+ \sum _{n<0} f_n(q^{-\frac{c}{2}}z)^{-n}\right),
\lb{Fminuss} 
\en here the positive (resp. negative) currents are expanded around $\infty$ (resp. 0).\
\noi
{\em {This yields the following decomposition of the total currents}}
\noi
\bea \lb{totaldecomp}
  &&\ \ a_1\eta (1)E(u)=E^+\left(u+\frac{c}{4}\right)-E^-\left(u-\frac{c}{4}\right),\nn\\
&&-a_2\eta(1)F(u)=F^+\left(u-\frac{c}{4}\right)-F^-\left(u+\frac{c}{4}\right).\qquad
\ena
\end{dfn}
 
\noi
One can derive the relations between these half-currents and the result is the following:

\begin{prop}[Commutation relations for Positive Half-Currents \cite{me1}]\lb{plusrelns}  
\bea
&& K_i^+(u_1)K_j^+(u_2)=K_j^+(u_2)K_i^+(u_1), \qquad   (i,j =1,2) \lb{plusrelkk} \\
&&\kpone^{-1}\Ept\kpone=\Ept\frac{q\eta(u-1)}{\eta(u)}+\Ep\frac{q^{-1}\eta(1)\eta(P+u-2)}{\eta(P-2)\eta(u)},
\lb{plusrelek1}\\
&&\kpone\fpt\kpone^{-1}=\fpt \frac{q\eta(u-1)}{\eta(u)} +\fp\frac{q^
{2u-1}\eta(1)\eta(P+h-u)}{\eta(P+h)\eta(u)},
\lb{plusrelfk1}\\
&& K_2^+ (u_1)^{-1 } E^+(u_2)K_2^+ (u_1)= E^+(u_2)\frac{q^{-1}\eta(u+1)}{{\eta }(u)}-E^+ (u_1)\frac{q^{-1}\eta   (1) \eta   (P+u)}{\eta   (P) \eta   (u)} , 
\lb{plusrelek2} \\
&&K_2^+(u_1)F^+(u_2)K_2^+(u_1){}^{-1}= F^+(u_2)\frac{q^{-1}\eta (u+1)}{\text{$\eta $}(u)} -\fp\frac{q^{2u-1}\eta (1)\eta (P+h-u-2)}{\eta (P+h-2)\eta (u)},\nn\\\lb{plusrelfk2}\\
&&\frac{q^{2u-1}\eta (1-u)}{\eta (u)}\hspace{0.2cm}E^+\left(u_1\right)E^+\left(u_2\right)+\frac{q^{-1}\eta (1+u)}{\text{$\eta $}(u)}\hspace{0.2cm}E^+\left(u_2\right)E^+\left(u_1\right) \nn\\
&&\hspace{0.5cm} =E^+\left(u_1\right){}^2\hspace{0.2cm}\frac{q^{-1}\eta (1)\eta (P+u-2)}{\text{$\eta $}(u)\eta (P-2)}+E^+\left(u_2\right){}^2\hspace{0.2cm}\frac{q^{2u-1}\eta (1)\eta (P-u-2)}{\text{$\eta $}(u)\eta (P-2)},\lb{plusrelee} \\
&&\frac{q^{-1}\eta (1+u)}{\eta (u)}\hspace{0.2cm} F^+\left(u_1\right) F^+\left(u_2\right)+\frac{q^{2u-1}\eta (1-u)}{\eta (u)}\hspace{0.2cm} F^+\left(u_2\right) F^+\left(u_1\right) \nn\\
&&\hspace{0.5cm}=F^+\left(u_1\right){}^2\hspace{0.2cm}\frac{q^{2u-1}\eta (1)\eta (P+h-u-2)}{\eta (P+h-2)\eta (u)}+ F^+\left(u_2\right){}^2\hspace{0.2cm}\frac{q^{-1}\eta (1)\eta (P+h+u-2)}{\text{$\eta $}(P+h-2)\eta (u)}, \\ \lb{plusrelff}
&&\left[E^+\left(u_1\right),F^+\left(u_2\right)\right] \nn\\
&&\quad=q^{2u}\left(q^{-1}-q\right)\left(K_2^+\left(u_2\right){}^{-1}K_1^+\left(u_2\right)\frac{\eta (P-u-1)}{\eta (u)\eta (P-1)} -K_2^+\left(u_1\right){}^{-1}K_1^+\left(u_1\right)\frac{\eta (P+h-u-1)}{\eta (u)\eta (P+h-1)}\right).\lb{plusrelef}\nn\\ 
\ena
\end{prop}
\subsection{The algebra \texorpdfstring{$U(R):=U(R_{q,x}(\slth))$}{U(R)}} Let us define the dynamical affine algebra based on the $RLL$ relations associated to a dynamical affine $R$-matrix.
The following $R$-matrix can be derived \cite{clad,idzumi} using the representation theory of $\U_q({\slth})$  and is the degeneration at $p \to 0$ of the elliptic $R$-matrix presented in \cite {JKOS2,K2}:

\bea \lb{rplus}
&&R_{q,x}^\pm(u,P)=\rho^\pm (z)\left(
\begin{array}{cccc}
 1 & 0 & 0 & 0 \\
 0 & \frac{\text{q$\eta $}(P-1)\eta (P+1)\eta (u)}{\eta (P)^2\eta (u+1)} & \frac{\eta (1)\eta (P+u)}{\eta (u+1)\eta (P)} & 0 \\
 0 & \frac{\eta (1)\eta (P-u)q^{2u}}{\eta (u+1)\eta (P)} & \frac{\text{q$\eta $}(u)}{\eta (u+1)} & 0 \\
 0 & 0 & 0 & 1 \\
\end{array}
\right). \nn\\
\ena

\noi
The expanded form is given by:

\bea
&&R_{q,x}^\pm(u,P)=\rho^\pm (z)\left(
\begin{array}{cccc}
 1 & 0 & 0 & 0 \\
 0 & \frac{q\left(1-q^{2(P-1)}\right) \left(1-q^{2(P+1)}\right)\left(1-q^{2u}\right)}{\left(1-q^{2P}\right)^2 \left(1-q^{2(u+1)}\right)} & \frac{\left(1-q^2\right) \left(1-q^{2(P+u)}\right)}{\left(1-q^{2P}\right) \left(1-q^{2(u+1)}\right)} & 0 \\
 0 & \frac{q^{2u}\left(1-q^2\right)\left(1-q^{2(P-u)}\right)}{\left(1-q^{2 P}\right)\left(1-q^{2(u+1)}\right)} & \frac{q\left(1-q^{2u}\right)}{1-q^{2(u+1)}} & 0 \\
 0 & 0 & 0 & 1 \nn\\
\end{array}
\right).
\ena
Here, the coefficient is defined to be:
\bea \lb{rho}
&&\rho^+ (z)=\frac{q^{\frac{1}{2}} \left(z^{-1};q^4\right)_\infty \left(q^4z^{-1};q^4\right)_\infty}{\left(q^2z^{-1};q^4\right)_\infty^2},\quad \text{}
 (a;x)_\infty=\prod_{k=0}^\infty(1-ax^{k})
\ena
and $\rho^-(z)=(\rho^+(z^{-1}))^{-1}$, which expands as
\bea \lb{rhominus}
&&\rho^- (z)=\frac{q^{-\frac{1}{2}}{\left(q^2 z;q^4\right)_\infty^2}}{ \left(z;q^4\right)_\infty \left(q^4z;q^4\right)_\infty}.
\ena
We will write $\rho^+(u)$ for $\rho^+(z)=\rho(q^{2u})$ when the context is clear.  We then have, $\rho^+ (u)\rho^- (-u)=1.$
The {\em elliptic R-matrix} which degenerates to $R_{q,x}(u,P)$ as $p \to 0$ is

\bea
&&R_{q,p}^+(u,P)=\widehat{\rho^+}(u)
\left(
\begin{array}{cccc}
1 &0                  &0             &0 \\
0  &\frac{\theta(P+1) \theta(P-1)\theta(u)  }{\theta(P)^2\theta(1+u)}    &\frac{\theta(1)\theta(P+u)}{\theta(P)\theta(1+u)} &0  \\
 0 &\frac{\theta(1)\theta(P-u)}{\theta(P)\theta(1+u)}&\frac{\theta(u)}{\theta(1+u)}&0 \\
0  & 0                 &0             &1 \\
\end{array}
\right),
\ena
\lb{Rmat4} where the Jacobi theta symbol is defined as
\bea \lb{thetadef}
&&\hspace{-1cm} \theta(u):=\frac{q^{\frac{u^2}{r}-u}}{(p;p)_\infty^3}\Theta_{p}(q^{2u}), {\text{ where}}\quad p=q^{2r}, \quad
\Theta_{p}(z)=(z;p)_\infty(p/z;p)_\infty(p;p)_\infty
\ena 
and the coefficient $\widehat{\rho^+}(u)$ becomes $\rho^+(u)$ as $p \to 0$.

\begin{dfn} [Extended Heisenberg algebra for $U(R)$]  \lb{heisenbergextended}
\bea
&&H=\C P + \C (P+h), \quad H_Q^*=\sum \C Q, \quad H^*=\frak{h}^* + H_Q^* \nn
\ena
\noi
We identify $\widehat{f}=f(P,P+h) \in\C [H]$ and  meromorphic functions on $\h^*$
via
\bea
&&f(P,P+h)(\xi)=f(<P,\xi>,<P+h,\xi>),\quad \xi \in H^*
\ena
\end{dfn}

\begin{dfn} \lb{UR}
 $U(R):=U(R_{q,x}(\glth))$ is the algebra over $M_{H^*}$ with generators $q^{\pm\frac{c}{2}}$ and $L^\pm(u)=\left(L^\pm_{ab}(u)\right)_{a,b=1}^2$, where $L^\pm_{ab}(u)=\sum_{n=0}^\infty L_{ab,\pm n}q^{\mp 2un}$,
and the RLL-relations:
\bea \lb{rll}
&&R^{\pm(12)}(u,P+h)L^{\pm (1)}\left(u_1\right)L^{\pm (2)}\left(u_2\right)=L^{\pm (2)}\left(u_2\right)L^{\pm (1)}\left(u_1\right)R^{\pm(12)}(u,P), \nn\\\
&&R^{\pm(12)}\left(u\pm\frac{c}{2},P+h\right)L^{\pm(1)}\left(u_1\right)L^{\mp(2)}\left(u_2\right)=L^{\mp(2)}\left(u_2\right)L^{\pm(1)}\left(u_1\right)R^{\pm(12)}\left(u\mp\frac{c}{2},P\right).\nn\\
\ena
\noi
\end{dfn}

\noi
{\em {Remarks.}} 
(i) Since the element $h$ is already in $U_q(\slth)$, Definition \ref{heisenbergextended} is consistent with the definition of $\mathbb{H}$ for  $U_{q,x}(\slth)$ in Section \ref{subsectionmaindef}.

\noi
(ii) $U(R)$ is an $H$-algebra with the $H$-bigrading defined by
\be
&& U(R)=\bigoplus_{\alpha,\beta \in H^*}{U(R)_{\alpha,\beta}}\hspace{.1cm},
\quad
 U(R)_{\alpha,\beta}=
\left\{\ x\in U(R) \left|\ \mmatrix{q^{P+h}xq^{-(P+h)}=q^{<\alpha, P>}x, \cr
q^{P}xq^{-P}=q^{<\beta, P>}x\cr}\ \right.\right\},
\en
and moment maps 
\bea
&&\mu_l(f)=f(P+h), \qquad   \mu_r(f)=f(P), \qquad f \in M_{H^*} 
\ena

\noi
(iii) The action on the generators
is given by
\bea
f(P+h)L^{\pm}_{ab}(u)&=&
L^{\pm}_{ab}(u)f(P+h-w(a)),\nn\\
f(P)L^{\pm}_{ab}(u)&=&L^{\pm}_{ab}(u)
f(P-w(b)),\lb{shifts}
\ena
where the weight function $w:\left\{1,2\right\} \rightarrow \left\{\pm1\right\}$ is given by identifying 
\bea
 L^{\pm}(u) \to \left(
\begin{array}{cc}
 L_{++}(u) & L_{+-}(u) \\
L_{-+}(u) & L_{--}(u) \\
\end{array}
\right)
\ena
\noi
(iv)
The dynamical tensor product of $U(R)$ with iteself, denoted $U(R) \tot U(R)$, is the $H^*$ bigraded vector space defined by
\be
 (U(R) {\widetilde{\otimes}}U(R))_{\al\beta}=\bigoplus_{\gamma\in\h^*} (U(R)_{\al\gamma}\otimes_{M_{\h^*}}U(R)_{\gamma\beta}),
\en
where $\otimes_{M_{\hbs}}$ is the usual tensor product $\otimes$ modulo the relation:
\bea \lb{dtpidentity}
&&\mu_r(f)a \tilde{\otimes } b = a \tilde{\otimes } \mu_l(f)b.
\ena
The Comultiplication $\Delta:U(R^\pm) \rightarrow U(R^\pm) \tot U(R^\pm)$ is given by  
\bea
&&\Delta \left(\mu _l(\hat{f})\right) = \mu _l(\hat{f})\tilde{\otimes }1,\qquad
\Delta \left(\mu _r(\hat{f}\right) = 1\tilde{\otimes }\mu _r(\hat{f}),\\
&&\Delta\left(e^Q\right)=e^Q \tot e^Q, \qquad \Delta \left(L^\pm_{\text{ab}}(u)\right)=\sum _{k=1}^2 L_{\text{ak}}^\pm(u)\tilde{\otimes }L_{\text{kb}}^\pm(u),
\ena
and the Counit $\varepsilon: U(R^\pm) \rightarrow \D$, the algebra of shift operators, is
\bea 
 &&\varepsilon \left(\mu _l(\hat{f})\right)=\varepsilon \left(\mu _r(\hat{f})\right)=\hat{f}T_0, \lb{epsilonUR1}\\
&&\varepsilon \left(e^Q\right)=e^Q, \qquad \varepsilon \left(L_{\text{ab}}^\pm(u)\right)=\delta _{a,b}T_{\text{w(b)Q}\text{  }}. \lb{epsilonUR2}
\ena
The expansion of the $RLL$-relations \eqref{rll} is given in Appendix A.
\begin{dfn}
The algebra $U(R_{q,x}(\slth))$ is the subalgebra of ${U(R)}$ with $Det(L^\pm(u))=1$, where the dynamical determinant element is given by 
\be 
&& Det(L^{\pm}(u))=L^{\pm}_{11}(u+1)L^{\pm}_{22}(u) - q^{}\frac{\eta(P-1)}{\eta(P)} L^{\pm}_{12}(u+1)L^{\pm}_{21}(u).
\en
\end{dfn}

\subsection{Equivalence of the two realizations}{\lb {equivdef}}
The following Ding-Frenkel type isomorphism (Theorem 4.4 in  \cite{me1}) will be the key to constructing dynamical representations of $U_{q,x}(\slth)$.  
\begin{thm} \lb{mainthm}
There exists a unique Gauss decomposition of the $L^\pm$-operators:
\bea \label{Lplusop}
&&L^\pm(u)=\left(
\begin{array}{cc}
 1 & F^\pm(u) \\
 0 & 1 \\
\end{array}
\right)\left(
\begin{array}{cc}
 K_1^\pm(u) & 0 \\
 0 & K_2^\pm(u) \\
\end{array}
\right)\left(
\begin{array}{cc}
 1 & 0 \\
 E^\pm(u) & 1 \\
\end{array}
\right) \\
&&\text{  }\text{  }\text{  }\text{  }\text{  } \text{  }\text{  } \text{  }\text{  }  =\left(
\begin{array}{cc}
 K_1^\pm(u)+F^\pm(u)K_2^\pm(u)E^\pm(u) & F^\pm(u)K_2^\pm(u) \\
 K_2^\pm(u)E^\pm(u) & K_2^\pm(u) \\
\end{array}
\right) \lb{gaussL}
\ena
which yields an an isomorphism of $H$-algebras
\be
&&\Phi : U_{q,x}\left(\slth\right) \longrightarrow U\left(R_{q,x}(\slth)\right), \text{ defined by}\\
&&\hspace{36pt} E(u) \mapsto  E^+\left(u+\frac{c}{4}\right)-E^-\left(u-\frac{c}{4}\right),\\
&&\hspace{36pt}F(u) \mapsto  F^+\left(u-\frac{c}{4}\right)-F^-\left(u+\frac{c}{4}\right),\\
&&\hspace{30pt}K_i^{\pm }(u) \mapsto  K_i^{\pm }(u), \quad q^{\pm c }\mapsto  q^{\pm c }.
\en
\end{thm}
Denote $U(R):=U(R_{q,x}(\slth)$.   The quantum loop algebra $\LL_{q}:=U_q(L(sl_2))$ can viewed as the quotient of the derived algebra of $U_{q}(\slth)$ by the 2-sided ideal generated by the element $c$.   Similarly define $\LL_{q,x}$ by replacing $U_q$ by $U_{q,x}$.

\begin{dfn}[Positive and Negative subalgebras]
$U(R^+)$ (resp. $U(R^-)$) is the subalgebra of $U(R)$ consisting of non-positive (resp. non-negative) powers in $z$.
Similarly, 
$\LL^+_{q,x}$ (resp. $\LL^+_{q,x}$) is the subalgebra of $\LL_{q,x}$ consisting of non-positive (resp. non-negative) powers in $z$.
\end{dfn}
\begin{prop}\lb{comult}
The coproduct for the half-currents is given as:
\bea
&&\Delta(K_1^{\pm}(u))=K_1^{\pm}(u)\tot K_1^{\pm}(u)+\sum_{j=1}^\infty (-1)^jK_1^{\pm}(u)E^{\pm}(u-1)^j\tot 
F^{\pm}(u-1)^jK_1^{\pm}(u),\nn\\
&&\Delta(K_2^{\pm}(u))=K_2^{\pm}(u) \tot K_2^{\pm}(u) + K_2^{\pm}(u)  E^{\pm}(u) \tot F^{\pm}(u) K_2^{\pm}(u),\nn\\
&&\Delta(E^{\pm}(u))=1\tot E^{\pm}(u)+E^{\pm}(u)\tot K_2^{\pm}(u)^{-1}K_1^{\pm}(u)\nn\\
&&\qquad\qquad\qquad\qquad+\sum_{j=1}^\infty(-1)^jE^{\pm}(u)^{j+1}\tot K_2^{\pm}(u)^{-1}F^{\pm}(u)^jK_1^{\pm}(u),\nn\\
\ena
\bea
&&\Delta(F^{\pm}(u))=F^{\pm}(u)\tot 1+ K_1^{\pm}(u)K_2^{\pm}(u)^{-1}\tot F^{\pm}(u)\nn\\\
&&\qquad\qquad\qquad\qquad+\sum_{j=1}^\infty(-1)^jK_1^{\pm}(u)E^{\pm}(u)^jK_2^{\pm}(u)^{-1}\tot F^{\pm}(u)^{j+1},\nn\\
&& \Delta \left(H^{+}(u)\right)=H^{+}(u)\hat{\otimes } H^{+}(u) \bmod A_{\geqslant 0}\tpt A_{\leqslant 0},
\ena 
where $A_{\geqslant 0}$  (resp. $A_{\leqslant0}$) is the subalgebra of $U(R)$ generated by the elements  $K_i^+(u),\ q^{\pm c}$ {\text{and }}
\noi
$ E^+(u)\ (\text{resp.} F^+(u)).$
\end{prop}

\section{Representation Theory of \texorpdfstring{$U_{q,x}(\slth)$}{Uqx}}

\subsection {Modules over (non-dynamical) quantum affine algebra \texorpdfstring{$U_q(\slth)$}{Uq} } 
\noi
Quantum affine $U_q=U_q(\slth)$ has a rich and  well-known representation theory \cite{CP}.   The main results, for our purposes, concerning the finite-dimensional representations are sumarized below.
Begin with the Poincare-Birkhoff-Witt formula, which reads, in standard notation
\bea \lb{pbw}
&&U_q= \mathcal{N}_-\mathcal{H}\mathcal{N}_+.
\ena
\begin{dfn}[Highest Weight Vector] \lb{qhwv}
A  highest weight vector $\Omega$ of a $U_q$-module $V$ is defined as a simultaneous eigenvector for all elements in $\mathcal{H}$ with
\bea \lb{qhw1}
&&e_{n }\Omega  = 0,\qquad \forall n\in \mathbb{Z}.
\ena
Such an $\Omega$ generates  a  {\em Highest Weight Module (HWM)}.  
\end{dfn}
\noi
{\em {Remarks}}
(i).
It is easy to verify (Proposition 3.2 in \cite{CP}) that any finite-dimensional module must be a HWM with $q^c=1$.
\noi
(ii)  Recall the quantum loop algebra $\LL_q$.   One may consider (see Prop 3.3 in \cite{CP}) only $\LL_q$-modules when constructing {\em finite-dimensional} $U_q(\slth)$-modules.

\begin{dfn} \lb{qhwvloop}
A highest weight vector $\Omega$ of highest weight $\underline{d} = \left\{d_n^{\pm }\right\}$ for an $\LL_q$-module is defined by
\bea\lb{qhw2}
&i)&e_k\Omega=0, \ \forall k\in \mathbb{Z},\nn\\
&ii)&\psi _n\Omega  = d_n^+ \Omega, \quad \varphi _{-n}\Omega  = d_{-n}^- \Omega , \text{  }\forall n\in \mathbb{N}\nn\\
&iii)&d_0^+d_0^- = 1.
\ena
\end{dfn}

From \eqref{pbw}, one obtains the induced (Verma) module $M(\underline{d})$ and its irreducible quotient  $V(\underline{d})$, in the standard way.   An important result we will need from \cite{CP} classifies the irreducible, finite-dimensional HWM's in terms of certain polynomials $P(t)$:
\begin{thm} [Chari-Pressley] \lb{qdpthm}
$V(\underline{d})$ is finite-dimensional $\iff  \exists  P(t)$, with 
$P(0) \neq 0$ and
\bea \lb{doubleqdp}
&&\sum _{n=0}^{\infty } d_n^+ t^n = q^{\deg  P}\frac{ P\left(q^{-2}t\right)}{P(t)} = \sum _{n=0}^{\infty } d_{-n}^-t^{-n},
\ena
where the left(right) equalities are expanded around 0($\infty$).
\end{thm}   
The polynomial $P(u)$ is known as the {\em {Drinfeld Polynomial}} associated to $V(\underline{d})$.
The main example of $V(\underline{d})$ is the Evaluation Module $V_l(w)$, based on the evaluation morphism from $U_q(\slth) \to U_q(sl_2)$ attributed to M. Jimbo \cite{Jimbo}.
\begin{dfn} [Evaluation Module $V_l(w)$ of $\LL_q$] \lb{qem}
The vectors $\left\{{v_i}\right\}_{i=0}^l$ form a weight basis for $V_l(w)$ with the generators acting as:
\bea
&&q^h v_i = q^{l-2i} v_i,\\
&&a_m v_i=\frac{w^m}{m}\frac{1}{q-q^{-1}}\left((q^m+q^{-m})q^{mh}-(q^{(l+1)m}+q^{-(l+1)m})\right)v_i,\lb{evaluationrep}\\
 &&e_n v_i = w^nq^{n(l-2i+1)}[l-i+1]{}_qv_{i-1}, \\
&&f_nv_i = w^nq^{n(l-2i-1)}[i+1]{}_qv_{i+1}, \qquad m,n \in \Z, m \neq 0. 
\ena
\end{dfn}
\noi
The corresponding Drinfeld polynomial is:
\bea
&&P(t)=\left(1-q^{n-1}{wt}\right)\left(1-q^{n-3}{wt}\right)\text{...}\left(1-q^{-n+1}{wt}\right).
\ena
{\em{Remark}.} It can be shown (Prop 4.3 in \cite{CP}) that the Drinfeld polynomials are well-behaved under tensor products:
whenever $V$ and $W$ are finite-dimensional representations of $\LL_q$ and their tensor product is irreducible,
we have
\bea
P_{V\otimes W}=P_VP_W,
\ena
Finally, the following result from \cite{CP} elegantly characterizes finite-dimensional $\LL_q$ modules as tensor products of the evaluation modules:
\begin{thm}[Chari-Pressley] \lb{qtpt}
For $i=1$ to $n$, there exist integers $l_i$ and  $w \in \mathbb{C} ^*$ such that
$V \in  \text{Rep}_f(\LL_q) \Longrightarrow  V \cong  V_{l_1}(w) \otimes V_{l_2}(w) \otimes  \text{...} V_{l_k}(w)$ uniquely, upto permuting the tensor factors.
\end{thm}
 
\subsection {The dynamical case \texorpdfstring{$U_{q,x}(\slth)$}{Uqx-mod}}
We will use the basic notions of dynamical representations of Hopf algebroids from \cite {KR,K2}.
Write $\C h$ as $\frak{h}_0$ and 
consider a vector space $\hV$ which is $\C h$-diagonalizable in the sense:
\be
&&\hV=\bigoplus_{\mu\in \frak{h}_0^*}\hV_{\mu},\qquad 
\hV_{\mu}=\{ v\in V\ |\ q^{x_0}v=q^\mu v\quad (x_0\in \frak{h}_0)\}.
\en

\noi
Define the $\h$-algebra $\cD_{\h,V}=\bigoplus_{\al,\beta\in \bH^*}(\cD_{\h,V})_{\al\beta}$, where
\be
&&\hspace*{-10mm}(\cD_{\h,V})_{\al\beta}=
\left\{\ X\in \End_{\C}\hV \left|\ \mmatrix{X(f(P)v)=f(P-<\beta,P>)X(v), \cr
X(V_\mu)\subseteq V_{\mu+\phi^{-1}(\al)-\phi^{-1}(\beta)}, f(P)\in \C[H] \cr}\ \right.\right\},
\en
along with moment maps
\be
&&\mu_l^{\cD_{H,V}}(\widehat{f})v=f(P+\mu)v,\quad 
\mu_r^{\cD_{H,V}}(\widehat{f})v=f(P)v\qquad 
\en
for $v\in V_{\mu}$.

\begin{dfn}[Dynamical representation]
A dynamical representation of $U_{q,x}'(\slth)$ on $\hV$ is 
 an $\h$-algebra homomorphism $\widehat{\pi}: U_{q,x}' 
 \to \cD_{\h,\hV}$. The dimension of the dynamical representation 
 $(\widehat{\pi},\hV)$ is  $\dim_{\frak{h}_0}\hV$.  
\end{dfn}
 
We will construct the dynamical representations of $U_{q,x}(\slth)$ in a similar fashion  to  the elliptic case {\cite {K2}} of $U_{q,p}(\slth)$.  In fact, most of our results for $U_{q,x}(\slth)$-modules coincide exactly with the corresponding results for $U_{q,p}(\slth)$ as $p \to 0$. Writing the $L^\pm(u)$ operators as
\be
L^\pm(u)=\left(
\begin{array}{cc}
L_{11}^{\pm} (u) &L_{12}^{\pm} (u) \\
 L_{21}^{\pm} (u) &L_{22}^{\pm} (u) \\
\end{array}
\right),
\en
the left action of  $P$ and $P+h$ is easily seen to be:
\begin{center}
    \begin{tabular}{ | l | l | l | l | l | l | l | l |}
    \hline
    $\text{ }$ & $L_{11}^{\pm} (u)$ & $L_{12}^{\pm} (u)$ & $L_{21}^{\pm}(u)$ & $L_{22}^{\pm} (u)$\\ \hline
    $P$ & $P-1$ & $P+1$ & $P-1$ & $P+1$\\ \hline
    $P+h$ & $P+h-1$ & $P+h-1$ & $P+h+1$ & $P+h+1$\\ \hline
     \end{tabular}
\end{center}
For example, $f(P)L_{11}^{\pm}(u)=L_{11}^{\pm}(u)f(P-1)$ and $f(P+h)L_{11}^{\pm}(u)=L_{11}^{\pm}(u)f(P+h-1)$.

\noi
By using a standard normal ordering procedure on the Heisenberg algebra, the Poincare-Birkhoff-Witt formula becomes
\bea \lb{pbwdynamical}
&&U'_{q,x}=({\cal N}_-\otimes{\cal H}\otimes{\cal N}_+)\otimes \C[\bH^*],
\ena
considered as a semi-direct product of $({\cal N}_-\otimes{\cal H}\otimes{\cal N}_+)$ and $\C[\bH^*]$.

 Consider now a vector space $V$ over  $\K=\C[H]$ which is $\C h$ diagonalizable.  Define $V_Q$ to be a vector space over $\C$ equipped with an action of $e^{ Q}$, and denote $\hV = V \otimes_{\Bbb{C}} V_Q$.
The vector space $\hV$ has an action of $f(P)$ and $e^Q$ via:
\be
&&f(P).(v\otimes \xi)=f(P)v\otimes \xi,\\
&&e^{Q}.(f(P)v\otimes \xi)=f(P+1)v\otimes e^{Q}\xi.
\en

\noi
We will only consider modules of this type hereafter.  Let us begin by defining the notion of pseudo-highest weight representations of $U_{q,x}((\slth))$ and of $U(R)$ and establishing their equivalence.

\begin{dfn} [Pseudo-highest weight representation of $U'_{q,x}(\slth)$]\lb{phwtotal}
A  dynamical representation $(\widehat{\pi}_V, \widehat{V}=V\otimes_{\Bbb{C}} V_Q)$ of  $U'_{q,x}(\slth)$ is said to be pseudo-highest weight, if there exists a vector (pseudo-highest weight vector)  $\hOmega$ in $\widehat{V}$,\text{ } $\hOmega= \Omega \otimes 1,  \text{  } \Omega \in V$ 
such that
\be
&i)&e^{ Q}.\hOmega=\hOmega\\
&ii)&e_k.\hOmega=0\quad \forall k \in \Bbb{Z}, \\
&iii)&\hOmega{\text{ is a simultaneous eigenvector for all the elements in }} \cal{H},\\
&iv)&\hV=U'_{q,x}(\slth).\hOmega.
\en
\end{dfn}

\begin{thm} \lb{totalhwmthm}
If $(\widehat{\pi}_V, \widehat{V}=V\otimes V_Q)$ is a finite-dimensional irreducible dynamical representation of $U'_{q,x}(\slth)$ then $\hV$ is a pseudo-highest weight representation of  $U'_{q,x}(\slth)$.  Further, $q^c$ acts as 1 or -1.
\end{thm}
\noi
{\it Proof.}
Consider the first statement of the theorem.  Since $V$ is finite-dimensional as a $\K[U'_q(\slth)]$-module, there is a vector $\Omega$ satisfying the conditions of Definition \ref{qhwv}, allowing us to define a vector $\hOmega':=\Omega \otimes \xi$ obeying conditions (ii) and (iii) in Definition \ref{phwtotal}.
 To prove condition (i) in Definition \ref{phwtotal}, observe that there are two types of $\Omega$:

\begin{itemize}
\item $\Omega$ is independent of P.
The action of  $e^{ Q}$  on $\hOmega'$ is 
$e^{ Q}.\hOmega'=\Omega\otimes e^{ Q}\xi$. Since $\hV$ is irreducible and finite-dimensional, there must exist a unique $\xi \neq 0$ and a complex constant $K \neq 0$ such that $e^Q\xi=K\xi$.  We can equate $\xi$ to $1$ by identifying $e^Q$ as $\frac{1}{K}e^Q$.

\item $\Omega$ depends on P.
Since $\hV$ is finite-dimensional, only a finite number of vectors in 
$\{\Omega(P+n)\ (n\in \Z)\}$ are $\K$-linearly independent.   By defining
$\widehat{\Omega}=\sum_{n\in \Z}\hOmega(P+n)\otimes \xi$, it is clear that
$e^{ Q}.\widehat{\Omega}=\sum_{n\in \Z}\hOmega(P+n)\otimes e^{ Q}\xi$. 
Then the same argument as the first case applies, so that $\xi=1$, yielding the required $\hOmega$ satisfying $(i)$. 
\end{itemize}
Condition (iv) is a consequence of \eqref{pbwdynamical}. 
Finally, the action of $q^c$ is proven in Corollary 3.2 in \cite{CP}.
\qed

Recall $\LL_{q,x}:=U_{q,x}(L(\slt)$ defined in Section \ref{equivdef}.
Observe that the $H$-Hopf algebroid structure on $U'_{q,x}(\slth)$ descends to $\LL_{q,x}$, since 
the bialgebroid structure and antipode are independent of $c$ by definition.
\begin{dfn} [Pseudo-highest weight representation of $\LL_{q,x}:=U'_{q,x}(L(\slt))$]\lb{phwloop}
A  dynamical representation $(\widehat{\pi}_V, \widehat{V}=V\otimes_{\Bbb{C}} V_Q)$ of $\LL_{q,x}$ is said to be pseudo-highest weight, if there exists a vector (pseudo-highest weight vector)  $\hOmega$ in $\widehat{V}$,\text{ } $\hOmega= \Omega \otimes 1,  \text{  } \Omega \in V$ and scalars 
  $\left\{d_n^{\pm }\right\}\in \C$ such that $d_0^+d_0^- = 1$, satisfying
\be
&i)&e^{ Q}.\hOmega=\hOmega,\\
&ii)&e_k.\hOmega=0\quad \forall k \in \Bbb{Z}, \\
&iii)&H^+_n(u).\hOmega= d_n^+ \hOmega, 
\quad H^-_{-n}(u).\hOmega= d_{-n}^- \hOmega,\qquad \forall n \in \Bbb{Z}_{\geq 0},\\
&iv)&\hV=\LL_{q,x}.\hOmega.
\en
\end{dfn}
\noi
This definition can be rephrased in terms of the algebra $U'(R)$, the quotient of $U(R)$ by the two-sided ideal generated by $q^{\pm\frac{c}{2}}$.
\begin{dfn}[Dynamical Pseudo-Highest Weight Modules of  $U'(R)$]\lb{phwrep}
A  dynamical representation $(\widehat{\pi}_V, \widehat{V}=V\otimes V_Q)$ of $U'(R)$ is said to be pseudo-highest weight, if there exists a vector (the pseudo-highest weight vector) $\hOmega\in \widehat{V}, \la\in \C$ and functions 
\bea
&&A^\pm(u)=\sum_{m \geq 0}A_{\pm m} z^{\mp m},
 D^\pm(u)=\sum_{m \geq 0}D_{\pm m} z^{\mp m},\ {\text {with }}\ A_{\pm m}, D_{\pm m} \in \C,\ {\text {such that }}\nn\\
&& A^\pm(u)=D^{\pm}(u-1)^{\pm1}  \nn
\ena
obeying the following conditions:
\be
&i)&e^{ Q}.\hOmega=\hOmega,\\
&ii)&L_{21}^{\pm}(u).\hOmega=0\quad \forall u, \\
&iii)&q^h.\hOmega=q^\la\hOmega,\quad L_{11}^{\pm}(u).\hOmega=A^{\pm}(u)\hOmega, 
\quad L_{22}^{\pm}(u).\hOmega=D^{\pm}(u)\hOmega, \\
&iv)&\hV=U'(R).\hOmega.
\en
\end{dfn}
 
\noi
We will denote the dynamical pseudo-highest weight as $(\lambda, A^\pm(u), D^\pm(u))$.    Then we have the following implications:
\begin{thm} \lb {dphwmthm}
\noi   Let  $\hV = V \otimes V_Q$. 
\begin{enumerate}
\item $(\widehat{\pi}_V, \widehat{V})$ is a DPHWM of $U'(R)$  if and only if it is a DPHWM over $\LL_{q,x}$.
\item If $\hV$ is a DPHWM of  $\LL_{q,x}$ then $V$ is a HWM over $\LL_q$.
\end{enumerate}
\end{thm}
 
\noi
{\it Proof of statement (i).}
Let $\hV$ be a DPHWM of $U'(R)$.  We must show that the conditions (ii) and (iii) in both definitions agree.  Condition (ii) is the same because of the series in Definition \ref{halfcurrentseries} for the half-currents.  Condition (iii) is the same because the Gauss decomposition \eqref{Lplusop} and the definitions of $H^\pm(u)$ in \eqref{hatefhdefn} yield
\bea
&&H^{\pm}(u). \Omega=A^{\pm}(u)A^{\pm}(u+1)\Omega.  
\ena
By expanding both sides as Laurent series in $z$,  one finds that the $\phi_m$ and $\psi_m$ are simultaneously diagonalized on $\Omega$ with the coefficients of the series on the right side determining their eigen values. 
The relation $d_0^+d_0^- = 1$ is equivalent to the observation that term
$ H^+_0H^-_0=q^hq^{-h}e^{2Q}=e^{2Q}$
acts as 1.

\noi
{\it Proof of statement(ii).}
All the conditions in Definition \ref{qhwvloop} are a part of Definition \ref{phwloop}.\qed 
\\

\noi
Now consider the induced (Verma) module $M(\vec{d})$ and its irreducible quotient  $V(\vec{d})$.
\begin{dfn}[Verma module]  
The Verma module $M({\vec d})$ is the quotient of  
$\LL_{q,x}$ by the left ideal generated by $\{x^+_k\ (k\in \Z),\ \psi_n-d^+_n\cdot 1, \phi_{-n}-d^-_{-n}\cdot 1\ (n\in \Z_{\geq 0})
,\  e^{ Q}-1
\}$.
\end{dfn}
\begin{prop}\lb{verma}
Given any sequence of complex-numbers $\vec{d}= \left\{d_{\pm n}^{\pm }\right\}_{n \ge 0}$,
the Verma module  $N(\vec{d})$ is a pseudo-highest weight representation of 
 pseudo-highest weight ${\vec d}$. Every pseudo-highest weight representation with pseudo-highest weight ${\vec d}$ is isomorphic to 
a quotient of $N(\vec{d})$. Furthermore, $N(\vec{d}))$ has a 
maximal proper submodule $M(\vec{d})$ such that
 $N(\vec{d})/ M(\vec{d})$
 is the unique irreducible pseudo-highest weight module of $\LL_{q,x}$, which is unique up to isomorphism. 
\end{prop}
We have the following dynamical extension of Theorem \ref{qdpthm}.
\begin{thm} \lb{ddpthm}
Consider the irreducible DPHWM $(\widehat{\pi}_V, \widehat{V}=V\otimes V_Q)$ of $\LL_{q,x}.$
A necessary and sufficient condition for $(\widehat{\pi}_V,\hV)$ to be finite-dimensional is:
$ \exists  P_V(u)$ (defined upto a scalar multiple), with 
$P_V(0)=1$ such that
\bea \lb{ddpcondition}
&&H^{\pm}(u).\hOmega=C_V\frac{ P_V\left(u+1\right)}{P(u)}.\hOmega,
\ena
where $C_V$ is a constant.
\end{thm} 
\noi  
{\it Proof.} We will prove the necessary part first and the sufficiency part after proving the Theorems \ref{rep} and \ref{TPDP}.  Recall Definition \ref{qhwvloop} and the remarks precdeding it.  Assuming finite-dimensionality of $\hV$, a pseudo highest-weight vector $\Omega$ for $V$ exists by the second statement in Theorem \ref{dphwmthm}.    By Theorem \ref{qdpthm}, there exists $P(z)$ such that the equation \eqref{doubleqdp} holds (with $z^{-1}$ replacing $t$ there) .  Factorizing $P(z)$ as $P(z)=\prod_{j=1}^{{\rm deg}P}(1-a_jz)$ and defining $P_V(u)=\prod_{j=1}^{{\rm deg}P}[u-\al_j]$,  where we denote $a_j=q^{2\al_j}$ and $[x]:=q^{-x}-q^x$, let us calculate the image of $H^\pm(u)$ under $\widehat{\pi}_{\hV}$ acting on $\hOmega$:
\be
H^\pm(u).\hOmega&=&q^{{\rm deg}P}\prod_{j=1}^{{\rm deg}P}
\frac{1-(a_j/q^2z)}{1-(a_j/z)}\hOmega\\
&=&q^{{\rm deg}P}\prod_{j=1}^{{\rm deg}P}\frac{q[u+1-\al_j]}{[u-\al_j]}\hOmega,
\en
as required.
\qed

\noi
We remark that the constant $C_V$ can be set equal to 1 by a gauge transformation {\footnote{eqn.(2.11) in \cite{JKOS2}} (as, for example, in Corollary \ref{ddpem}). We will refer to $P(u)$ as the {\it Dynamical Drinfeld Polynomial (DDP).}

\subsection{Dynamical Evaluation Modules}

The main example of a finite-dimensional dynamical highest weight representation is the dynamical evaluation module.  The action of the currents can be explicitly derived, similar to the non-dynamical and elliptic cases.
We will use the following notations henceforth:
\bea
&&[x]_q:=\frac{q^{x}-q^{-x}}{q-q^{-1}}, \quad  [x]:={q^{-x}-q^{x}},  \quad{\text {so that}}\quad \frac{[x]_q}{[y]_q}=\frac{[x]}{[y]}.
\ena
 
Thus, $\lim_{p \rightarrow 0}{\theta(u)}=[u]$ and the connection with Gasper-Rahman is
$\lim_{p \rightarrow 0}{[u]_{GR}}=[u]_{q^{1/2}}$ (see eqns. (1.6.9) and (11.2.3) in \cite {GR})

For $a \in \C$, the evaluation representation $V_l(q^{2a})$ of $L_q$ given in Definition \ref{qem} becomes a dynamical representation by setting
\be
\widehat{V}:=V^{(l)}(q^{2a}) = V_l(q^{2a}) \otimes V_Q, \quad {\text {where }} V_Q=\C1, \ {\text {and  }} e^Q(f(P)v \otimes 1)=f(P+1)v \otimes 1.
\en
The next result explicitly describes the Dynamical Evaluation Module structure and can be viewed as the dynamical extension of Definition \ref{qem} from the previous section.

\begin{thm} {\bf Dynamical Evaluation Module of  $\LL_{q,x}$ and  $\LL^{\pm}_{q,x}$. } \lb{demtotal}
 Let $V^{(l)}(q^{2a})=V_l(q^{2a}) \otimes \C 1$ and define the operators $S^{\pm}$ by $S^{\pm}(v_m) = v_{m \mp 1}$.  Then  $(\widehat{\pi}_{l},V^{(l)}(q^{2a}))$ is a dynamical representation of the (sub)algebras
$\LL_{q,x}, \LL^{+}_{q,x}$ and $\LL^{-}_{q,x}$ via the following actions:

(i) $\LL_{q,x}$:
\bea
\widehat{\pi}_{l,w}(E(u))&=&S^+\frac{[\frac{l+h+2}{2}]}{{[1]}}\delta \left(u-a-\frac{h+1}{2}\right)e^{2Q},\\
\widehat{\pi}_{l,w}(F(u))&=&S^-\frac{[\frac{l-h+2}{2}]}{{[1]}}\delta \left(u-a-\frac{h-1}{2}\right),\\
\widehat{\pi}_{l,w}(H^{\pm}(u))&=&\frac{[u-a-\frac{l+1}{2}][u-a+\frac{l+1}{2}]}{{[u-a-\frac{h-1}{2}][u-a-\frac{h+1}{2}]}}e^{2Q},
\ena

(ii)  $\LL^{+}_{q,x}$: 
\bea
\widehat{\pi}_{l,w}(K_1^{+} (u)) &=& \rho^+_l(u-a-1)^{-1}\frac{[u-a+\frac{l-1}{2}]}{[u-a-\frac{h+1}{2}]}e^Q,\\
\widehat{\pi}_{l,w}(K_2^{+} (u)) &=& \rho^+_l(u-a)\frac{[u-a-\frac{h-1}{2}]}{[u-a+\frac{l+1}{2}]}e^{-Q},\\
\widehat{\pi}_{l,w}(E^{+}(u))&=&-e^Q S^+\frac{[u-a-\frac{h+1}{2}-P][\frac{l+h+2}{2}]}{{[u-a-\frac{h+1}{2}][P]}}e^{Q}, \\
\widehat{\pi}_{l,w}(F^{+}(u))&=&S^-\frac{[u-a+\frac{h-1}{2}+P][\frac{l-h+2}{2}]}{{[u-a-\frac{h-1}{2}][P+h-1]}},
\ena

(iii) $\LL^{-}_{q,x}$:
\bea
\widehat{\pi}_{l,w}(K_1^{-} (u)) &=& \rho^-_l(u-a-1)^{-1}\frac{[u-a+\frac{l-1}{2}]}{[u-a-\frac{h+1}{2}]}e^Q,\\
\widehat{\pi}_{l,w}(K_2^{-} (u)) &=& \rho^-_l(u-a)\frac{[u-a-\frac{h-1}{2}]}{[u-a+\frac{l+1}{2}]}e^{-Q},\\
\widehat{\pi}_{l,w}(E^{-}(u))&=&e^Q S^+\frac{[-u+a+\frac{h+1}{2}+P][\frac{l+h+2}{2}]}{{[-u+a+\frac{h+1}{2}][-P]}}e^{Q}, \lb{repeminus}\\
\widehat{\pi}_{l,w}(F^{-}(u))&=&-S^-\frac{[-u+a-\frac{h-1}{2}-P][\frac{l-h+2}{2}]}{{[-u+a+\frac{h-1}{2}][-P-h+1]}},\lb{repfminus}
\ena
where 
\be
&&\rho_l^+ (z)=\frac{q^\frac{l}{2} \left(q^{l+3}z^{-1};q^4\right)_\infty \left(q^{-l+1}z^{-1};q^4\right)_\infty}{\left(q^{3-l}z^{-1};q^4\right)_\infty \left(q^{l+1}z^{-1};q^4\right)_\infty}, \quad
\rho_l^- (z)=\frac{q^{-\frac{l}{2}}\left(q^{3-l}z^{};q^4\right)_\infty \left(q^{l+1}z^{};q^4\right)_\infty }{\left(q^{l+3}z^{};q^4\right)_\infty \left(q^{-l+1}z^{};q^4\right)_\infty},
\en
so that $\rho^-_l(z)^{-1}=\rho^+_l(z^{-1})$.
\end{thm}
The formulae \eqref{repeminus} and \eqref{repfminus} should be understood as expansions around $z=0$ where $z=q^{2u}$.   Notice that the decompositions of the total currents into half-currents \eqref{totaldecomp} remains valid at the module level, by using \eqref{etadelta}. 

{\em Proof}. 
 The action of the dynamical Drinfeld currents and half-currents is derived using the expressions in Definitions \ref{hatefhdefn}, \ref{halfcurrentseries} and \ref{halfcurrentintegrals} in the (non-dynamical) evaluation module in Definition \ref{qem}, recalling that the minus (plus) half-currents are realized as series around 0 ($\infty$).
One can directly verify the commutation relations in Propositions \ref{Defrelns} and  \ref{plusrelns} for the action on $\widehat{V}:=V^{(l)}(a) = V_l(a) \otimes V_Q$ by using the identity
\bea \lb {etaid3}
&&\hspace{-.5cm}\frac{\eta(u_1+t)\eta(u_2 + s)}{\eta(u_1) \eta(u_2) \eta(s)} =  \frac{\eta(u_1- u_2+t)\eta(u_2 + s + t)}{\eta(u_1 - u_2) \eta(u_2) \eta(s + t)}+
\frac{\eta(u_2- u_1 + s) \eta(u_1 + s + t) \eta(t)}{\eta(u_2-u_1)\eta(u_1) \eta(s) \eta(s + t)}.\nn\\
\ena
 
\qed

\noi
The module $(\widehat{\pi}_{l},V^{(l)}(q^{2a}))$ will be referred to as the {\it Dynamical 
Evaluation Module (DEM)}.  We are pleased to observe that the formulae are the same as for the elliptic quantum group $U_{q,p}(\slth)$, as $p \rightarrow 0$, in Theorem (4.13) in \cite{K2} and equations (C.8) and (C.9) in \cite{JKOS2}.  The next result describes the action of the matrix elements of the $L^+(u)$-operator on this vector space.
\begin{thm} {\bf Dynamical Evaluation Module of  $U(R^+)$} \lb{rep}
The dynamical action of the matrix elements of $L^+(u)$ on  $V^{(l)}(a) = V_l(a) \otimes V_Q$ is given by
\be
\widehat{\pi}_{l,w}(L^+_{11}(u))&=&-\frac{[u-a+\frac{h+1}{2}][P-\frac{l-h}{2}][P+\frac{l+h+2}{2}]}{\varphi^+_l(u-a)[P][P+h+1]}e^Q,\\
\widehat{\pi}_{l,w}(L^+_{12}(u))&=&-S^-\frac{[u-a+\frac{h-1}{2}+P][\frac{l-h+2}{2}]}{\varphi^+_l(u-a)[P+h-1]}e^{-Q},\\
\widehat{\pi}_{l,w}(L^+_{21}(u))&=&S^+\frac{[u-a-\frac{h+1}{2}-P][\frac{l+h+2}{2}]}{\varphi^+_l(u-a)[P]}e^Q,\\
\widehat{\pi}_{l,w}(L^+_{22}(u))&=&-\frac{[u-a-\frac{h-1}{2}]}{\varphi^+_l(u-a)}e^{-Q},
\en
where
\be
&&\varphi^+_l(u)=-\rho_{l}^+(u)^{-1}[u+\frac{l+1}{2}].
\en
\end{thm}
\noindent {\it Proof}.
It is straightforward to derive these relations using the half-current integrals in Definition\ref{halfcurrentintegrals} and Theorem \ref {demtotal}. The Riemann addition formula \eqref{raf} is used for $\widehat{\pi}_{l,w}(L^+_{11}(u))$.   In the course of verifying the $RLL$-relations, we use the formula
\bea \lb{phirel}
\varphi^+_l(u) \varphi^+_l(u-1)=[u-\frac{l+1}{2}][u+\frac{l+1}{2}],
\ena
which can be checked by using the expression $\rho_l^+(q^{2u})$ for $\rho_l^+(u)$.
\qed \\

\noi
{\it Remark.  }
$V^{(l)}(a)$ is a DPHWM over $U(R^+)$ with pseudo-highest weight
\be
&&(\lambda, A^+(u), D^+(u))=\left( l,-\frac{[u-a+\frac{l+1}{2}]}{\varphi^+_l(u-a)}, -\frac{[u-a-\frac{l-1}{2}]}{\varphi^+_l(u-a)} \right).
\en
\begin{cor} \lb{repminus}
The vector space $\hV(a)$ is a $U(R^-)$-module by
\be
\widehat{\pi}_{l,w}(L^-_{11}(u))&=&\frac{[a-u+\frac{h+1}{2}][P-\frac{l-h}{2}][P+\frac{l+h+2}{2}]}{\varphi^-_l(u-a)[P][P+h+1]}e^Q,\\
\widehat{\pi}_{l,w}(L^-_{12}(u))&=&S^-\frac{[a-u-\frac{h-1}{2}-P][\frac{l-h+2}{2}]}{\varphi^-_l(u-a)[P+h-1]}e^{-Q},\\
\widehat{\pi}_{l,w}(L^-_{21}(u))&=&-S^+\frac{[a-u+\frac{h+1}{2}+P][\frac{l+h+2}{2}]}{\varphi^-_l(u-a)[P]}e^Q,\\
\widehat{\pi}_{l,w}(L^-_{22}(u))&=&\frac{[a-u+\frac{h-1}{2}]}{\varphi^-_l(u-a)}e^{-Q},
\en
where 
\be
&&\varphi^-_l(u)=-\rho_{l}^-(u)^{-1}[u+\frac{l+1}{2}].
\en
\end{cor}
\noi
{\em Proof}.  Straightforward, since \eqref{phirel} remains true on replacing $\varphi_l^+(u)$ by $\varphi_l^-(u)$.  The pseudo-highest weight is  $(\lambda, A^-(u), D^-(u))$=$\left( l,\frac{[a-u-\frac{l+1}{2}]}{\varphi^-_l(u-a)}, \frac{[a-u+\frac{l-1}{2}]}{\varphi^-_l(u-a)} \right)$.
\qed

\begin{cor} \lb{ddpem}
$\hV(a)$ is a DPHWM over $U(R)$ with psuedo-highest weight $(\lambda, A^\pm(u), D^\pm(u))$.
The corresponding dynamical Drinfeld polynomial is
\bea
&&P(u)=[u-a-\frac{l-1}{2}][u-a-\frac{l-1}{2}+1]\text{ ... }[u-a+\frac{l-1}{2}].
\ena
\end{cor}
\noi
{\em Proof.}  The second set of $RLL$-relations \eqref{rll}, expanded in Appendix (A.1) can be verified - since the central element $c$ acts as zero, they are essentially the same as the first set.    It is straightforward to check that $P(u)$ satisfies the required condition \eqref{ddpcondition}.
\qed

Since we get a dynamical representation of the {\em total} algebra as well as each of the {\em half-current} subalgebras (i.e. the dynamical affine quantum groups) on the  same underlying vector space, without loss of generality, we will hereafter consider the representations of $\LL_{q,x}^{+}$ only. 

\section{Tensor Products}
\subsection{Construction for the Evaluation Modules}
Recall the dynamical tensor product on $U(R)$ given in \eqref{dtpidentity}:
\bea \lb{dtprel}
&&f(u,P)a \tilde{\otimes } b = a \tilde{\otimes } f(u,P+h)b.
\ena
The tensor product of dynamical representations becomes a dynamical representation due to the important observation that there exists a natural $\h$-algebra embedding from $\theta_{VW}: 
\cD_{\h,\hV}\tot \cD_{\h,\hW}\to \cD_{\h,\hV\tot \hW}$, where the dyanamical tensor product is given by
\be
&&\hV\tot \hW=\bigoplus_{\al\in \bar{\hh}^*}(\hV\tot \hW)_{\al},\quad (\hV \tot \hW)_{\al}=\bigoplus_{\beta\in  \bar{\hh}^*}\hV_{\beta}\otimes_{M_{\h^*}}\hW_{\al-\beta}.
\en
Here
$\otimes_{M_{H^*}}$ denotes the usual tensor product modulo the relation
\bea
&&f(P)v\otimes w=v\otimes f(P+\nu)w\lb{VtotW}
\ena
for $w\in \hW_\nu$.  Finally, the tensor product $\hV\tot \hW$ acquires an action of $\C[\hbs]$ through:
\be
&&f(P).(v\tot w)=\Delta(\mu_r(\widehat{f}))(v\tot w)=v\tot f(P)w.
\en

The next result, that the DDP continues to be well-behaved under tensor products, just like in the ordinary quantum situation, will be used in the proof of the sufficiency part of Theorem \ref{ddpthm}.
\begin{thm} \lb {TPDP}
Let  $(\widehat{\pi}_V,\hV)$ and $(\widehat{\pi}_W,\hW)$ be finite-dimensional dynamical modules of $\LL_{q,x}$ such that $\hV \tot \hW$ is irreducible.  Then, 
\bea
&& P_{\hV\tot \hW}=P_{\hV}P_{\hW}.
\ena
\end{thm}
\noi
{\em Proof.}  The comultiplication formulae for the half-currents in Proposition \ref{comult} are used to verify this result.\qed

We can now prove the {\em sufficiency part of Theorem \ref{ddpthm}.}
Let $P_{\hV}(u)$ be any function satisfying the conditions of the theorem, factored as
$P_{\hV}(u)=\prod_{j=1}^r{[u-a_j]}.$ We will define $\hV$ as follows.  By Theorem \ref{dphwmthm}, $P_{\hV}(u)$ uniquely determines $\vec{d}$, the set of eigenvalues of
$\left\{ \psi_n,\varphi_{-n} \right\} _{n \ge 0}$.
For each $a_j$, let  $\tilde{V}_j=V^{(1)}(q^{2\alpha_j})$ denote the $U(R)$-module from Corollary \ref{repminus}, with PHWV $\tilde{v}_j:=v_0^{(1)}\otimes{1}$ and let  $W=\otimes_j{\tilde{V}_j}$.
Evidently, ${W}$ is a DPHWM with PHWV $\Omega:=\otimes_j{\tilde{v}_j}$, satsifying the condition that $q^h\Omega=q^r\Omega$.   Now, the $\LL_{q,x}$-submodule $W'$ of $W$ generated by $\Omega$ has a maximal submodule $W''$ and the quotient module $W'/W''$ is irreducible and by Corollary \ref{ddpem} and Theorem \ref{TPDP}, it has DDP $P_{\hV}(u)=\prod P_j(u)$, defined upto a scalar multiple.  This allows us to define $\hV$ as $W'/W''$.   It is clearly a finite-dimensional irreducible DPHWM as required.
\qed

The following result, a degeneration of the elliptic one (Prop 4.16 in \cite{K2}), confirms that the DEM of Theorem \ref{rep} coincides with the $R-$matrix representation. 
\begin{prop}\lb{LandR}
Let us define the matrix elements of $\widehat{\pi}_{l,w}(\hL^\pm_{\vep_1\vep_2}(u))$ by
\be
\widehat{\pi}_{l,w}(\hL^\pm_{\vep_1\vep_2}(u))v^l_m&=&
\sum_{m'=0}^l (\hL^\pm_{\vep_1\vep_2}(u))_{\mu_{m'}\mu_m}v^l_{m'},
\en
where $\mu_m=l-2m$. Then we have
\be
(\hL^\pm_{\vep_1\vep_2}(u))_{\mu_{m'}\mu_{m}}=
R^+_{1l}(u-v,P)_{\vep_1 \mu_{m'} }^{\vep_2 \mu_m}.
\en
Here $R^+_{1l}(u-a,P)$ is the degeneration of the $R$-matrix from (C.17)  in  \cite{JKOS2}, as $p \to 0$. 
The  case $l=1$, $R^+_{11}(u-a,P)$ coincides with the image 
$(\pi_{1,z}\otimes \pi_{1,w})$ of the universal $R$-matrix $\cR^+(u,P)$\cite{JKOS2} given in \eqref{rplus}.
The case $l>1$, $R^+_{1l}(u-a,P)$ coincides with the $R$-matrix 
obtained by fusing $R^+_{11}(u-a,P)$ $l$-times.  For $\widehat{\pi}_{l,w}(\hL^-_{\vep_1\vep_2}(u))$, just replace $\rho^+(u-a)$ by $\rho^-(u-a)$  for $R^+_{1l}(u-a,P)$ in the previous calculation to obtain  $R^-_{1l}(u-a,P)$ . 
\end{prop}

Now we can establish a dynamical version of the final result from the previous subsection, characterising the tensor product modules.  An elliptic version of this result appears in \cite{K2} and \cite{FV}.  Denote the shifted $q$-factorial as
$[u]_m=[u][u+1]\cdots[u-m+1]$.

\begin{thm} \lb{dtpt}
The tensor product $V:=V^{(l_1)}(a) \otimes V^{(l_2)}(b)$ is a DPHWM over $U(R^+)$ if and only if $b-a=\frac{l_1+l_2-2s}{2}+1$, \hspace{0.1cm} $ s=0,1\cdots\rm{min}(l_1,l_2)$.   Explicitly, the pseudo-highest weight $(\la, A(u), D(u))$ and pseudo-highest vector $\Omega^{(s)}_V$ are given by
\be
&&\Omega^{(s)}_V =\sum_{m=0}^{s}{C_{m}^{s}(P) v_m^{l_1} \tilde{\otimes}  v_{s-m}^{l_2}}, \qquad
\la=l_1+l_2-2s, \\
&& C_{m}^{s}(P)= C_{0}^{s}\frac{[P-l_2+s-2m]_s [l_2-s+1]_m[P-2m+1]_{2m}}{[-l_1]_m [P-l_2+s-2m]_m [P+s-2m+1]_m}, \\
&& A(u) = \frac{[u-a+\frac{l_1+1}{2}][u-b+\frac{l_2+1}{2}+s]}{\phi_{l_1}(u-a) \phi_{l_2}(u-b)}   ,\qquad
D(u)=  \frac{[u-a-\frac{l_1-1}{2}+s][u-b-\frac{l_2-1}{2}]}{\phi_{l_1}(u-a) \phi_{l_2}(u-b)},
\en
where $C_0^s$ is a constant factor, independent of $P$.  
\end{thm}
\noi
{\it Proof.}  The verification of the theorem is along the same lines as the proof of Theorem 4.17 in \cite{K2} for the elliptic case.       The first condition in  
 Definition \ref{phwrep} follows from the formula $\Delta(e^Q)=e^Q \otimes e^Q$ and the construction of the dynamical representations $V^{(l_1)}(a)$ in Theorem \ref{rep}.
For the specified $\la$,  condition (iii) gives the decomposition of the pseudo highest-weight vector as
\bea \lb{omegasv}
\Omega^{(s)}_V&=&\sum_{m_1=0}^{min\{ l_1,s\}\ }C^s_{m_1}(u,P)v^{l_1}_{m_1}\tot v^{l_2}_{s-m_1}.
\ena
The coefficients can be determined using the conditions (ii) and (iv) by solving a recurrence relation.  More specifically, solving condition (ii) yields the following recurrence relation on the constants $C^s_{m_1}(u,P)$, which is independent of $u$ if and only if  $b-a=\frac{l_1+l_2-2s}{2}+1$:
\bea
C^s_{m_1}(u,P)&=&-C^s_{m_1-1}(u,P-2)\frac{[u-a-\frac{l_1+1}{2}+m_1][u-b
+\frac{l_2-2s+1}{2}+m_1+1-P]}{[u-a-\frac{l_1+1}{2}+m_1+1-P][u-b
+\frac{l_2-2s+1}{2}+m_1]}\nn\\
&&\quad \times\frac{[l_2-s+m_1][P][P-1]}{[l_1+1-m_1][P-l_2+s-1-m_1][P+s-m_1]}\nn.\lb{recursionc}
\ena
 The condition (iv) for $\delta(u)$ is used to reduce this to the recurrence relation:
\be
C_{m_1}^s(P)=C^s_{m_1}(P-1)\frac{[P][P-l_2+2s-2m_1-1]}{
[P-l_2+s-m_1-1][P+s-m_1]}.
\en
In the process one also obtains the formula for $D(u)$, by employing the well-known Riemann addition formula  
\bea \lb{raf}
[u'+x][u'-x][v+y][v-y]-[u'+y][u'-y][v+x][v-x]=[x-y][x+y][u'+v][u'-v],\nn\\
\ena
with $u'=u-a-\frac{l_1+l_2-2s}{2},$  $x = \frac{l_2+1}{2},$  $v=\frac{2P+2s-l_2-2m_1-1}{2}$ and $y=\frac{l_2-2s+2m_1+1}{2}$.

The calculation of $A(u)$ is similar.  Finally, the relation $A(u)D(u-1)=1$ can be verified by using \eqref{phirel}.
\qed\\

We can finally state the main result, confirming a conjecture of Konno \cite{K2}, that the intertwiner of tensor products of representations of the elliptic quantum group $U_{q,p}(\slth)$ (given by an elliptic ${}_{12}V_{11}$ series) degenerates precisely to a ${}_{10}W_9$ sum on the basis of the representation theory of the quantum group $U_{q,x}(\slth)$. 
\subsection{Hypergeometric Series and Intertwiners}
 Let us recall some basics regarding hypergeometric functions from the standard text {\cite{GR}}, using the base $q^2$ instead of $q$ and series in $z=q^2$.  \noi
We use the notations
\bea
&&(a;q^2)_k=\prod_{n=0}^{k-1}(1-aq^{2n}), \quad
(a_1,\cdots,a_{r+1};q^2)_k : = (a_1;q^2)_k \cdots (a_{r+1};q^2)_k.
\ena
\begin{dfn}\lb{basichyp} {The Basic Hypergeometric series $ {}_{r+1}\phi_r$ is defined as}
\bea \lb{hypphidef}
 {}_{r+1}\phi_r
\left(\mmatrix{a_1, \cdots,a_{r+1}\cr
b_1,\cdots,b_r\cr};q^2,q^2\right)=\sum_{k=0}^{\infty}\frac{(a_1,\cdots,a_{r+1};q^2)_k}{(q^2,b_1,\cdots,b_r;q^2)_k} \ q^{2k}.
\ena
\end{dfn}
\noi

\noi
(i) This series, ${}_{r+1}\phi_r$ is {\em well-poised} if $q^2a_1=a_2b_1=a_3b_2=\cdots =a_{r+1}b_r$. 

\noi
(ii) A well-poised series is {\em very well-poised} if $a_2=q^2 a_1^{\frac{1}{2}}$ and $a_3=-q^2a^{\frac{1}{2}}$. 

\noi
(iii) The {\em very well-poised balancing condition}  is given by  (eq(2.1.12) in \cite{GR}) 
\be 
&&(b_1 \cdots b_{r-2})q^2=(\pm a_1^{\frac{1}{2}}q)^{r-3}.
\en
\begin{dfn}
Wilson's \cite{Wilson} biorthogonal symbol ${}_{r+1}W_r$ is defined via
\bea \lb{hypWdef}
&& \quad{}_{r+1}W_r(a;b_1,\cdots,b_{r-2};q^2,q^2) =  {}_{r+1}\phi_r
\left(\mmatrix{a,q^2 a^{\frac{1}{2}},- q^2 a^{\frac{1}{2}},b_1,\cdots,b_{r-2}\cr
a^{\frac{1}{2}},- a^{\frac{1}{2}},{aq^2}/{b_1},\cdots,{aq^2}/{b_{r-2}}\cr};q^2,q^2\right)\nn\\
&&\qquad \qquad \qquad  \qquad \qquad \qquad \quad  =\sum_{k=0}^{\infty}\left(\frac{1-aq^{4k}}{1-a}\right) \frac{(a,b_1,\cdots,b_{r-2};q^2)_k}{(q^2,aq^2/b_1,\cdots,aq^2/b_{r-2};q^2)_k} \ q^{2k}.\nn\\
\ena
\end{dfn}

\begin{dfn} Recall the Jacobi theta functions $\Theta_p(u)$ and $\theta(u)$ in \eqref{thetadef}.  Define a very well-poised theta hypergeometric series ${}_{r+1}V_r$ by 

(i) Multiplicative form

\bea
\qquad\qquad{}_{r+1}V_r(a;b_1,\cdots,b_{r-4};q^2,p) =\sum_{k=0}^{\infty}\left(\frac{\Theta_{p}(aq^{4k})}{\Theta_{p}(a)}\right) \frac{(a,b_1,\cdots,b_{r-4};q^2)_k}{(q^2,aq^2/b_1,\cdots,aq^2/b_{r-4};q^2)_k} \ q^{2k}.\nn\\
\ena

(ii) Additive form
\bea
{}_{r+1}v_r(u_0;u_1,\cdots,u_{r-4})=\sum_{k=0}^\infty\frac{\theta(u_0+2k)}{\theta(u_0)}
\prod_{i=0}^{r-4}\frac{\theta(u_i)_k}{\theta(u_0+1-u_i)_k},
\ena
with the balancing condition
\be
\sum_{i=1}^{r-4}u_i=\frac{r-7}{2}+\frac{r-5}{2}u_0,
\en
which guarantees that the two forms are equal (see eq(11.3.8) and eq(11.3.25) in \cite{GR}).
\end{dfn}
\noi
It is well known that the non-elliptic degeneration of ${}_{r+1}V_r$ involves a shift by 2 in r:
\bea
&&\lim_{p \to 0}  {}_{r+1}V_r(a;b_1,\cdots,b_{r-4};q^2,p)= {}_{r-1}W_{r-2}(a;b_1,\cdots,b_{r-4};q^2,q^2)
\ena
\noi
Finally, {\em Jackson's summation formula} \eqref{jack} is true if the LHS is {\em VWP balanced}, i.e. if $a^2q^{2(k+1)}=bcde$.  
\bea \lb{jack}
 {}_{8}\phi_7
\left(\mmatrix{a,q^2a^{\frac{1}{2}},-q^2a^{\frac{1}{2}},b,c,d,e,q^{-2k}\cr
a^{\frac{1}{2}},-a^{\frac{1}{2}},aq^2/b,aq^2/c,aq^2/d,aq^2/e,aq^{2(k+1)}\cr};q^2,q^2\right)=\frac{(aq^2,aq^2/bc,aq^2/bd,aq^2/cd;q^2)_k}{(aq^2/b,aq^2/c,aq^2/d,aq^2/bcd;q^2)_k}.\nn\\
\ena

Let us prove our final result that  
the expression of a general basis element is governed by Wilson's  ${}_{10}W_9$ series as conjectured by Konno \cite{K2}:
\begin{thm}\lb{thm10w9}
If $l=l_1+l_2-2s$, we have, for $0\leq m\leq l$,
\bea
&&\Delta( L^+_{12}(u) L^+_{12}(u+1)\dots L^+_{12}(u+m-1))\ \Omega^{(s)}_V\nn\\
&&=\frac{[P]}{\prod_{i=1}^m\varphi_{l_1}(u-a+i-1)\varphi_{l_2}(u-a-
\frac{l}{2}+i-1)}\nn\\
&&\times\hspace*{-10mm}
\sum_{k={\rm max}(0,s+m-l_2)}^{{\rm min}(l_1,s+m)}
\hspace*{-6mm}(-1)^kC_0^s\frac{[P+m-2k-l_2+s]_s}{[P+m-2k+1]_s}
\frac{[u-a+\frac{l_1+1}{2}]_{m-k}[P+m-2k]}
{[u-a-l+\frac{l_1-1}{2}+m-k+P]_{m-k}[P-k]}
\nn\\
\ena
\bea
&&\qquad\times
\frac{[-u+a-\frac{l_1-1}{2}-m+k-P]_{k}[-u+a+l-\frac{l_1-1}{2}-m+1]_k}{[P+m-k][P+l_1-2k+1]_m} \nn\\
&&\qquad\times\frac{[-m]_k
[P-k]_{m-k}[P+l_1-k+1]_{m-k}[s+1]_{m-k}}
{[P-k+1]_{m-k}[P]_{m-k}}\nn\\
&&\qquad\times 
{}_{10}W_{9}\left(q^{2(P+m-2k)};q^{-2s},q^{-2k},q^{2(P-k)},q^{2(l_2-s+1)},q^{2(-u+a-\frac{l_1-1}{2})},\right.\nn\\
&&\qquad\qquad \qquad \left. q^{2(u-a-l+
\frac{l_1-1}{2}+2m-2k+P)},q^{2(P+m-2k+l_1+1)};q^2,q^2\right) \ v^{l_1}_k \tot\ v^{l_2}_{m+s-k},\lb{betavs}
\ena
while $\Delta( L^+_{12}(u) L^+_{12}(u+1)\dots L^+_{12}(u+m-1)) \ \Omega^{(s)}_V = 0$ for $m>l$.
\end{thm}
\noi
{\em Proof.}  We can work in the additive setting because of the commutativity of the following diagrams:

\begin{tikzpicture}
  \matrix (m) [matrix of math nodes,row sep=4em,column sep=6em,minimum width=2em] {
     {}_{12}V_{11} &{}_{12}v_{11} \\
     {}_{10}W_{9} &{}_{10}w_{9}\\
     };
  \path[-stealth]

    (m-1-1) edge node [left] {$p \to 0$} (m-2-1)
            edge [<->] node [above] {$VWP$} node [below] {=} (m-1-2)
         
    (m-2-1.east|-m-2-2) edge  [<->] node [below] {$VWP$} node [above] {=} (m-2-2)
                                
    (m-1-2) edge node [right] {$p \to 0$} (m-2-2)
;
\hspace{8cm}

  \matrix (m) [matrix of math nodes,row sep=3em,column sep=4em,minimum width=2em] {
    \Theta_p(x) & \theta(x) \\
     \eta(x) &\left[x\right]\\
     };
  \path[-stealth]

    (m-1-1) edge node [left] {$p \to 0$} (m-2-1)
            edge [<->] node [above] {} node [below] {} (m-1-2)
         
    (m-2-1.east|-m-2-2) edge  [<->] node [below] {} node [above] {} (m-2-2)
                                
    (m-1-2) edge node [right] {$p \to 0$} (m-2-2)
;
\end{tikzpicture}

The very well-poised condition (iii) in Definition \ref{basichyp} guarantees the commutativity of the above diagram (here $[x]:=q^{-x}-q^{x}$). 
  
Claim that the element of $U(R^+) \tot U(R^+)$ acting in the LHS of the theorem can be expanded as
\bea \lb{Deltabeta}
&&\Delta (L^+_{12}(u)L^+_{12}(u+1)\cdots L^+_{12}(u+m-1))=\sum_{j=0}^m B^m_j(P) X_{j,m}(u)\tot Y_{j,m}(u), \ {\text {where}} \nn\\
&&\qquad B^m_j(P)=\frac{[1]_m}{[1]_j[1]_{m-j}}\frac{[P][P-m+2j]}{[P+j][P-m+j]},\qquad  m\in \Z_{\geq 0},\nn\\
&&\qquad X_{j,m}(u)=L^+_{11} (u+m-1)\cdots L^+_{11} (u+m-j)L^+_{12}(u+m-j-1)\cdots L^+_{12}(u),\nn\\
&&\qquad Y_{j,m}(u)=L^+_{22}(u)\cdots L^+_{22}(u+m-j-1)
L^+_{12}(u+m-j)\cdots L^+_{12}(u+m-1).\nn\\
\ena
To prove the claim, we use induction on the variable $m$.  Let $m=1$ then 
\be
B_0^1(P)=B_1^{1}(P)=1, \
X_{0,1}(u)=L^+_{11}(u), \ X_{1,1}(u)=L^+_{12}(u), \quad  Y_{0,1}(u)=L^+_{12}(u), \ Y_{1,1}(u)=L^+_{22}(u),
\en
yielding
\be
\Delta(L_{12}^+(u))=L_{11}^+(u) \tot L_{12}^+(u) + L_{12}^+(u) \tot L_{22}^+(u),
\en
which is exactly the comultiplication formula for $L_{12}^+(u)$. 

 Let us assume the result holds for $m$ and prove the induction step by deriving a recurrence relation for $B_j^m(P)$. 

 By definition,
\be
&&\Delta (L^+_{12}(u)L^+_{12}(u+1)\cdots L^+_{12}(u+m))\nn\\
&&=\Delta (L^+_{12}(u)L^+_{12}(u+1)\cdots L^+_{12}(u+m-1))\Delta( L^+_{12}(u+m)) \nn\\
&&=\Delta (L^+_{12}(u)L^+_{12}(u+1)\cdots L^+_{12}(u+m-1))\left(L_{11}^+(u) \tot L_{12}^+(u) + L_{12}^+(u) \tot L_{22}^+(u)\right)\nn\\
&&=(\sum_{j=0}^m B^m_j(P) X_{j,m}(u)\tot Y_{j,m}(u))\left(L_{11}^+(u) \tot L_{12}^+(u) + L_{12}^+(u) \tot L_{22}^+(u)\right)\nn\\
&&=\left(\sum_{j=0}^m B^m_j(P) X_{j,m}(u) L_{11}^+(u) \tot Y_{j,m}(u)L_{12}^+(u) \right)
+ \left(\sum_{j=0}^m B^m_j(P) X_{j,m}(u) L_{12}^+(u) \tot Y_{j,m}(u)L_{22}^+(u) \right)\nn\\
\en
\be
&&=L^+_{11}(u+m-1)\cdots L^+_{11}(u)L^+_{11}(u+m) \tot L^+_{12}(u)\cdots L^+_{12}(u+m-1)L^+_{12}(u+m) \nn\\
&&+L^+_{12}(u+m-1)\cdots L^+_{12}(u)L^+_{12}(u+m) \tot L^+_{22}(u)\cdots L^+_{22}(u+m-1)L^+_{22}(u+m)\\
&&+\sum_{j=1}^m \left\{{\ \atop\ } D^m_{j-1}(P)
L^+_{11}(u+m-1)\cdots L^+_{11}(u+m-j+1)L^+_{12}(u+m-j)\cdots L^+_{12}(u)L^+_{11}(u+m)\right.\\
&&\left.+\frac{[P+j-m]}{[P+2j-m]}B^m_{j}(P)L^+_{11}(u+m-1)\cdots L^+_{11}(u+m-j)L^+_{12}(u+m-j-1)
\cdots
L^+_{12}(u)L^+_{12}(u+m)\right\}\\
&&\tot L^+_{22}(u)\cdots L^+_{22}(u+m-j-1)L^+_{22}(u+m-j)L^+_{12}(u+m-j+1)\cdots L^+_{22}(u+m).
\en
In the last step, we changed $j$ to $j-1$ in the first summand of the preceding equation and used the property \eqref{dtprel} along with the formula 
\be
L_{22}^+(u)L_{12}^+(u+1)=\frac{[P+h+1]}{[P+h]}L_{12}^+(u)L_{22}^+(u+1)
\en
obtained from (A.10).  The first two summands in the final line correspond to $j=m$ and $j=0$, in the first and second sum in the preceding line, respectively.
Compare with 
\be
\sum_{j=0}^{m+1} B^{m+1}_j(P) X_{j,m+1}(u)\tot Y_{j,m+1}(u).
\en
We need to prove that 
\bea \lb{needalphabeta}
&&B^{m+1}_j(P-j+1)L^+_{11}(u+m)L^+_{12}(u+m-j)\cdots L^+_{12}(u)\nn\\
&&=B^m_{j-1}(P-j+1)L^+_{12}(u+m-j)\cdots L^+_{12}(u)L^+_{11}(u+m)\nn\\
&&+\frac{[P-m+1]}{[P+j-m+1]}B^m_{j}(P-j+1)L^+_{11}(u+m-j)L^+_{12}(u+m-j-1)\cdots L^+_{12}(u)L^+_{12}(u+m) \nn\\
\ena
for $j=1,2,\cdots,m$.

\noi
Applying (A.3) successively yields the formula:
\bea \lb{alphabetaidentity}
L^+_{11}(u)L^+_{12}(v_1)\cdots L^+_{12}(v_l)
&=&\frac{[P+1][P-l][u-v_l]}{[P][P-l+1][u-v_1+1]}L^+_{12}(v_1)\cdots L^+_{12}(v_l)L^+_{11}(u)\nn\\
&&+\sum_{k=1}^l\frac{[P+1][P-k+1-u+v_k][1]}{[P][u-v_1+1][P-k+2]}L^+_{12}(v_1)\cdots L^+_{11}(v_k)\cdots L^+_{12}(v_l)L^+_{12}(u).\nn\\
\ena
Use \eqref{alphabetaidentity}, with $l\to m-j, u\to u+m-j, v_k\to u+k-1\ (1\leq k\leq m-j)$, to simplify the second term in the RHS of \eqref{needalphabeta}.  Use \eqref{alphabetaidentity} again with  $l\to m-j+1, u\to u+m, v_k\to u+k-1\ (1\leq k\leq m-j+1)$ in the LHS of \eqref{needalphabeta}.  Comparing coefficients produces the required recurrence relation
\be
\frac{B^m_{j}(P)}{B^{m}_{j-1}(P)}=\frac{[P-m+j-1][P-m+2j][m-j+1][P+j-1]}
{[P-m+j][P-m+2(j-1)][P+j][j]}.
\en
Solving this relation yields the formula for $B^m_{j}(P)$ as claimed.


Theorem \ref{rep} yields the action of $X_{j,m}(u)$ and $Y_{j,m}(u)$ on the basis elements. Let $k=m_1+m-j$. Then, 
\bea
&&\hspace{-1cm}X_{j,m}(u)v^{l_1}_{m_1}=(-1)^{k+m_1+m}\frac{[u-a+\frac{l_1+1}{2}]_{m-k}[P-k]_{m-k}
[P+l_1-k+1]_{m-k}[-u+a-m-\frac{l_1-1}{2}-P+k]_{k}}
{\prod_{i=1}^m\varphi_l(u-a+i-1)\ [P]_{m-k}[P+l_1-2k+1]_{m}}\nn\\
&&\qquad\qquad \times 
\frac{[1]_k[-u+a-\frac{l_1+1}{2}+1]_{m_1}[P-2k+m]_{m_1}[P+l_1+m-2k+1]_{m_1}}
{[P+m-k]_{m_1}[u-a+m+\frac{l_1-1}{2}+P-2k+1]_{m_1}[1]_{m_1}}\ v^{l_1}_{k},
\ena
and
\bea
&&v^{l_1}_{k}\tot Y_{j,m}(u)\ v^{l_2}_{s-m_1}\nn\\
&&=(-1)^{m+k}\frac{[-u+b+\frac{l_2-1}{2}-m-s+1]_{k}
[u-b-\frac{l_2+1}{2}+m+s+1-k+P]_{m-k}[s+1]_{m-k}}
{\prod_{i=1}^m\varphi_{l_2}(u-b+i-1)[P-k+1]_{m-k}}\nn\\
&&\quad\times \frac{[u-b-\frac{l_2+1}{2}+2m+s+1-2k+P]_{m_1}[P-k+1]_{m_1}[-s]_{m_1}}
{[-u+b+\frac{l_2-1}{2}-m-s+1]_{m_1}[P+m-2k+1]_{2m_1}}
v^{l_1}_k\ \tot\  v^{l_2}_{m+s-k}.
\ena
The first expression is obtained using Theorem \ref{rep} and noting that, for $b, k, m_1 ,k-m_1\in \Z_{\geq 0}$, 
\be
&&[a]_{m_1+b}=[a]_{b}[a+b]_{m_1},\qquad 
\\&&
[a-b]_{m_1+b}=(-1)^b[-a+1]_{b}[a]_{m_1},    \qquad \\
&&[a+k]_{k-m_1}=(-1)^{m_1}\frac{[a+k]_k}{[-a-2k+1]_{m_1}},\hspace*{7cm}
\en
while the second follows by employing \eqref{dtprel} and
\be
&&[a-k]_{k-m_1}=(-1)^{k+m_1}\frac{[-a+1]_k}{[-a+1]_{m_1}},\\
&&[a+m_1-k]_{m_1+m-k}=\frac{[a-k]_{m-k}[a+m-2k]_{2m_1}}{[a-k]_{m_1}}\qquad k,m,m_1\in 
\Z_{\geq 0},\ m\geq k\geq m_1.
\en

\noi
We can now calculate the action of the elements in \eqref{Deltabeta} on the pseudo-highest weight vector $\Omega_V^{(s)}$ \eqref{omegasv}.

\be
&& \Delta( L^+_{12}(u) L^+_{12}(u+1)\dots L^+_{12}(u+m-1))\Omega^{(s)}_V \nn\\
&&\qquad =\sum_{m_1=0}^s\sum_{j=0}^m
B_j^m(P)C_{m_1}^{s}(P-m+2j) X_{j,m}(u)
v^{l_1}_{m_1} \tot Y_{j,m}(u) v^{l_2}_{s-m_1}.
\en

Change the summation from $j$ to $k=m+m_1-j$ and use the following identity to expand the coefficient $C_{m_1}^s(P-m+2j)$ 
\be
&&\frac{[P-m+2j-l_2+s-m_1]_{s-m_1}}{[P-m+2j+1]_{s-m_1}}\nn\\
&&=\frac{[P-m+2j-l_2+s-2m_1]_s[P-m+2j-2m_1+1]_{2m_1}}{[P-m+2j-2m_1+1]_s[P-m+2j-l_2+s-2m_1]_{m_1}[P-m+2j+s-2m_1+1]_{m_1}}.
\en
Then, summation over $m_1$, for ${\rm max}(0,k-m)\leq m_1\leq {\rm min}(k,s)$, yields the ${}_{10}W_{9}$ part of the required result.  The remaining coefficients are the same.

For the second statement of the theorem, let $m=l+1$.  Using $k=l_1-s+n$, we have
\bea 
&&{}_{10}W_{9}\left(P+m-2k;-k,-s,P-k,l_2-s+1,-u+a-\frac{l_1-1}{2},\right.\nn\\
&&\qquad\qquad \qquad \qquad \left. u-a-l+
\frac{l_1-1}{2}+2m-2k+P,P+m-2k+l_1+1\right)\\
&&=\sum_{m_1=0}^s\frac{[P-l_1+l_2+1-2n+2m_1][P-l_1+l_2+1-2n]_{m_1}[-l_1+s-n]_{m_1}[-s]_{m_1}}{[P-l_1+l_2+1-2n][l_2-s+2-n]_{m_1}[P-l_1+s+1-2n]_{m_1}[-l_1]_{m_1}}\nn\\
&&\times\frac{[P-l_1+s-n]_{m_1}[l_2+1-s]_{m_1}[P+l_2+2-2n]_{m_1}}{[1]_{m_1}[P+l_2-s+2-n]_{{m_1}}[P-l_1+l_2+2+s-2n]_{m_1}}.\lb{419}
\ena

\noi
The VWP balancing condition is easily checked:
\be 2a+1 = 2P - 2 l_1 + 2 l_1 - 4 n + 3=b+c+d+e.
\en
Now apply the Jackson formula \eqref{jack}, in additive form, 
with the identifications:
\be
x \to q^{2x}, \quad  k \to s, \quad a=P-l_1+l_2+1-2n, \quad b=-l_1+s-n, \\
c=P-l_1+s-n, \quad  d=l_2+1-s, \quad e=P+l_2+2-2n,
\en
 so that \eqref{419} becomes 
\be
\frac{[P-l_1+l_2-2n+2]_s[l_1+l_2-2s+2]_s[P-n+1]_s[1-n]_s}{[P+l_2-n-s+2]_s[l_2-n-s+2]_s[P-l_1-2 n+P+s+1]_s[l_1-s+1]_s}
\en
which clearly vanishes for $n=1,2,\cdots,s$ since the last term on the numerator
\be
[1-n]_s=0  \quad  n=1,2,\cdots,s.
\en 
This completes the proof of the theorem.
\qed

\vspace{-5pt}
\begin{figure}[here]
\hspace{2cm}
\includegraphics[width=0.7\textwidth]{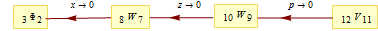}
\caption{Degeneration of elliptic 6j symbols}
\label{fig:hypfig}
\end{figure} 
\subsection{Concluding Remarks and Discussion}

The ${}_{10}W_{9}$ symbol is obtained on the basis of corepresentations (equivalently, representations) of elliptic $U(n)$ by Koelink et. al. in \cite{KNR}, so it seems natural that a version of the quantum Kazhdan-Lusztig functor (Figure \ref{fig:figure1.1}) exists between $U_{q,x}(\slth)$-mod and elliptic $U(2)$-mod.  The series in the first two boxes on the left in Figure \ref{fig:hypfig} have been derived in \cite{KR, KirRes, KoeKoo, Koorn} and correspond to the usual $q$-6j symbols and $q$-Racah coefficients.

The phenomenon of self-duality\cite{Rosengren} can viewed as an equivalence between modules and comodules. It is {\em dynamical} in nature, since the usual quantum algebras $F_q[SL_2]$ and $F_q[\widehat{SL_2}]$ have {\em finite-dimensional irreducible representations only of dimension 1} \cite{Soibelman,meold}, unlike the dynamical case.  The $F_{q,\la}[SL_2]$ modules in \cite{KR} and the current article exist in any dimension.  Thus $U(R)$ cannot contain the Kac-Moody dual of $U_q(\slth)$. 
 
We have established here the representation theory for finite-dimensional modules of $U_{q,x}(\slth)$, corresponding to level-zero modules of $U_q(\slth)$.
The infinite-dimensional ones have also been derived for the elliptic quantum group $U_{q,p}(\slth)$ by Konno et al. using dynamical $Z$-algebras and quantum $W$-algebras \cite{Knew}.  They exist at $p=0$ also and the precise relation with $U_{q,x}(\slth)$ will be presented elsewhere.

It would be interesting to find a categorical framework for our results, similar to the work of Shibukawa and Takeuchi \cite{ST} where an explicit isomorphism of tensor categories $A(R)$-mod $\simeq$ $R$-mod {\em as sets} is proven,  for $R=R(\la)$ in the finite (trigonometric) case.

\subsection*{Acknowledgments}
The author would like to thank George Andrews, Pavel Etingof, Hitoshi Konno, Nicolai Reshetikhin, Hjalmar Rosengren and Ping Xu for valuable collaborations and support. This work was completed at MSRI, Berkeley and at Pennsylvania State University.   

\begin{appendix}
\section {The $RLL$ relations of $U(R)$}
\begin{prop} \lb{rllexpanded}  Using the notation in \eqref{rappendix}, the first relation in \eqref{rll}
\be
&&R^{\pm(12)}(u,P+h)L^{\pm (1)}\left(u_1\right)L^{\pm(2)}\left(u_2\right)=L^{\pm (2)}\left(u_2\right)L^{\pm (1)}\left(u_1\right)R^{\pm(12)}(u,P),
\en
 is expanded as
\bea
&&L^{\pm}_{11} \left(u_1\right) L^{\pm}_{11} \left(u_2\right)=L^{\pm}_{11} \left(u_2\right) L^{\pm}_{11} \left(u_1\right),\quad
L^{\pm}_{12} \left(u_1\right) L^{\pm}_{12} \left(u_2\right)=L^{\pm}_{12} \left(u_2\right) L^{\pm}_{12} \left(u_1\right),\\
&&L^{\pm}_{21} \left(u_1\right) L^{\pm}_{21} \left(u_2\right)=L^{\pm}_{21} \left(u_2\right) L^{\pm}_{21} \left(u_1\right) ,\quad
L^{\pm}_{22} \left(u_1\right) L^{\pm}_{22} \left(u_2\right)=L^{\pm}_{22} \left(u_2\right) L^{\pm}_{22} \left(u_1\right),\\
&&L^{\pm}_{11} \left(u_1\right) L^{\pm}_{12} \left(u_2\right)=L^{\pm}_{11} \left(u_2\right) L^{\pm}_{12} \left(u_1\right) \bar{c}(u,P)+L^{\pm}_{12} \left(u_2\right) L^{\pm}_{11} \left(u_1\right) b(u,P),\lb{alphabeta}\\ 
&&L^{\pm}_{12} \left(u_1\right) L^{\pm}_{11} \left(u_2\right)=L^{\pm}_{11} \left(u_2\right) L^{\pm}_{12} \left(u_1\right) \bar{b}(u,P)+L^{\pm}_{12} \left(u_2\right) L^{\pm}_{11} \left(u_1\right) c(u,P),\\
&&L^{\pm}_{21} \left(u_1\right) L^{\pm}_{22} \left(u_2\right)=L^{\pm}_{21} \left(u_2\right) L^{\pm}_{22} \left(u_1\right) \bar{c}(u,P)+L^{\pm}_{22} \left(u_2\right) L^{\pm}_{21} \left(u_1\right) b(u,P),\\
&&L^{\pm}_{22} \left(u_1\right) L^{\pm}_{21} \left(u_2\right)=L^{\pm}_{21} \left(u_2\right) L^{\pm}_{22} \left(u_1\right) \bar{b}(u,P)+L^{\pm}_{22} \left(u_2\right) L^{\pm}_{21} \left(u_1\right) c(u,P), \\
&&b(u,P+h) L^{\pm}_{11} \left(u_1\right) L^{\pm}_{21} \left(u_2\right)+ c(u,P+h) L^{\pm}_{21} \left(u_1\right) L^{\pm}_{11} \left(u_2\right)=L^{\pm}_{21} \left(u_2\right) L^{\pm}_{11} \left(u_1\right), \nn\\ \\
&&b(u,P+h) L^{\pm}_{12} \left(u_1\right) L^{\pm}_{22} \left(u_2\right)+ c(u,P+h) L^{\pm}_{22} \left(u_1\right) L^{\pm}_{12} \left(u_2\right)=L^{\pm}_{22} \left(u_2\right) L^{\pm}_{12} \left(u_1\right), \nn\\ \\
&&\bar{b}(u,P+h) L^{\pm}_{21} \left(u_1\right) L^{\pm}_{11} \left(u_2\right)+\text{  }\bar{c}(u,P+h) L^{\pm}_{11} \left(u_1\right) L^{\pm}_{21} \left(u_2\right)=L^{\pm}_{11} \left(u_2\right) L^{\pm}_{21} \left(u_1\right), \nn\\ \\
&&\bar{b}(u,P+h) L^{\pm}_{22} \left(u_1\right) L^{\pm}_{12} \left(u_2\right)+\text{  }\bar{c}(u,P+h) L^{\pm}_{12} \left(u_1\right) L^{\pm}_{22} \left(u_2\right)=L^{\pm}_{12} \left(u_2\right) L^{\pm}_{22} \left(u_1\right), \nn\\ \\
&&b(u,P+h) L^{\pm}_{11} \left(u_1\right) L^{\pm}_{22} \left(u_2\right)+c(u,P+h) L^{\pm}_{21} \left(u_1\right) L^{\pm}_{12} \left(u_2\right) \nn\\
&&\qquad \qquad =L^{\pm}_{21} \left(u_2\right) L^{\pm}_{12} \left(u_1\right) \bar{c}(u,P)+L^{\pm}_{22} \left(u_2\right) L^{\pm}_{11} \left(u_1\right) b(u,P),\\
&&b(u,P+h) L^{\pm}_{12} \left(u_1\right) L^{\pm}_{21} \left(u_2\right)+c(u,P+h) L^{\pm}_{22} \left(u_1\right) L^{\pm}_{11} \left(u_2\right) \nn\\
&&\qquad \qquad =L^{\pm}_{21} \left(u_2\right) L^{\pm}_{12} \left(u_1\right) \bar{b}(u,P)+L^{\pm}_{22} \left(u_2\right) L^{\pm}_{11} \left(u_1\right) c(u,P),\\
&&\bar{b}(u,P+h) L^{\pm}_{21} \left(u_1\right) L^{\pm}_{12} \left(u_2\right)+\bar{c}(u,P+h) L^{\pm}_{11} \left(u_1\right) L^{\pm}_{22} \left(u_2\right) \nn\\
&&\qquad \qquad =L^{\pm}_{11} \left(u_2\right) L^{\pm}_{22} \left(u_1\right) \bar{c}(u,P)+L^{\pm}_{12} \left(u_2\right) L^{\pm}_{21} \left(u_1\right) b(u,P),\lb{A.13}\\
&&\bar{b}(u,P+h) L^{\pm}_{22} \left(u_1\right) L^{\pm}_{11} \left(u_2\right)+\bar{c}(u,P+h) L^{\pm}_{12} \left(u_1\right) L^{\pm}_{21} \left(u_2\right) \nn\\
&&\qquad \qquad =L^{\pm}_{11} \left(u_2\right) L^{\pm}_{22} \left(u_1\right) \bar{b}(u,P)+L^{\pm}_{12} \left(u_2\right) L^{\pm}_{21} \left(u_1\right) c(u,P).
\ena
The second set :
\be
&&R^{\pm(12)}\left(u\pm\frac{c}{2},P+h\right)L^{\pm(1)}\left(u_1\right)L^{\mp(2)}\left(u_2\right)=L^{\mp(2)}\left(u_2\right)L^{\pm(1)}\left(u_1\right)R^{\pm(12)}\left(u\mp\frac{c}{2},P\right),
\en

\vspace{-5pt}
\noi
is expanded below (we write the $R$-matrix entries in \eqref{rappendix}, as $c(u,P) = c_0(u,P)$ and  $\bar{c}(u,P) = \bar{c_0}(u,P)$  to distinguish them from the central element $c)$.

\bea
&& \rho^{\pm}(u \pm \frac{c}{2}) \ L^{\pm}_{11}(u_1)  L^{\mp}_{11}(u_2)= L^{\mp}_{11}(u_2)  L^{\pm}_{11}(u_1) \ \rho^{\pm}(u \mp \frac{c}{2}) ,\\
&& \rho^{\pm}(u \pm \frac{c}{2}) L^{\pm}_{11}(u_1)  L^{\mp}_{12}(u_2) \nn\\
&&\qquad \quad  \qquad  =\Big( L^{\mp}_{11}(u_2)  L^{\pm}_{12}(u_1) \bar{c_0}(u \mp \frac{c}{2},P)+ \  L^{\mp}_{12}(u_2) L^{\pm}_{11}(u_1) b(u \mp \frac{c}{2},P)\Big) \rho^{\pm}(u \mp \frac{c}{2}),\nn\\ \\
&& \rho^{\pm}(u \pm \frac{c}{2})  L^{\pm}_{12}(u_1)  L^{\mp}_{11}(u_2) \nn\\
&& \qquad \quad \qquad  = \Big(L^{\mp}_{11}(u_2)  L^{\pm}_{12}(u_1) \bar{b}(u \mp \frac{c}{2},P)+ \  L^{\mp}_{12}(u_2) L^{\pm}_{11}(u_1) c_0(u \mp \frac{c}{2},P) \Big)\rho^{\pm}(u \mp \frac{c}{2}),\nn\\ \\
&&\rho^{\pm}(u \pm \frac{c}{2})  L^{\pm}_{12}(u_1)  L^{\mp}_{12}(u_2)= L^{\mp}_{12}(u_2)  L^{\pm}_{12}(u_1) \rho^{\pm}(u \mp \frac{c}{2}),\\
&&\rho^{\pm}(u \pm \frac{c}{2})\Big(b(u \pm \frac{c}{2},P+h) L^{\pm}_{11}(u_1)  L^{\mp}_{21}(u_2) + c_0(u \pm \frac{c}{2},P+h)  L^{\pm}_{21}(u_1)  L^{\mp}_{11}(u_2)\Big)\nn\\  
&&\qquad \quad  \qquad =L^{\mp}_{21}(u_2)  L^{\pm}_{11}(u_1)  \rho^{\pm}(u \mp \frac{c}{2}),\\
&&\rho^{\pm}(u \pm \frac{c}{2}) \Big( b(u \pm \frac{c}{2},P+h) L^{\pm}_{11}(u_1) L^{\mp}_{22}(u_2)+ c_0(u \pm \frac{c}{2},P+h)  L^{\pm}_{21}(u_1)  L^{\mp}_{12}(u_2)\Big)\nn\\
&&\qquad \quad \qquad   = \Big(L^{\mp}_{22}(u_2) L^{\pm}_{11}(u_1) b(u \mp \frac{c}{2},P)+ L^{\mp}_{21}(u_2)  L^{\pm}_{12}(u_1) \bar{c_0}(u \mp \frac{c}{2},P) \Big)\rho^{\pm}(u \mp \frac{c}{2}),\nn\\ \\
&&\rho^{\pm}(u \pm \frac{c}{2}) \Big( b(u \pm \frac{c}{2},P+h)   L^{\pm}_{12}(u_1)  L^{\mp}_{21}(u_2)+ c_0(u \pm \frac{c}{2},P+h) L^{\pm}_{22}(u_1)  L^{\mp}_{11}(u_2)\Big)\nn\\
&&\qquad \quad  \qquad  = \Big(L^{\mp}_{21}(u_2)  L^{\pm}_{11}(u_1) {c_0}(u \mp \frac{c}{2},P) + L^{\mp}_{21}(u_2)  L^{\pm}_{12}(u_1) \bar{b}(u \mp \frac{c}{2},P)  \Big)\rho^{\pm}(u \mp \frac{c}{2}),\nn\\ \\
&&\rho^{\pm}(u \pm \frac{c}{2})\Big( b(u \pm \frac{c}{2},P+h)  L^{\pm}_{12}(u_1)  L^{\mp}_{22}(u_2)+c_0(u \pm \frac{c}{2},P+h)  L^{\pm}_{22}(u_1)  L^{\mp}_{12}(u_2)\Big)\nn\\
&&\qquad \quad  \qquad = L^{\mp}_{22}(u_2)  L^{\pm}_{12}(u_1) \rho^{\pm}(u \mp \frac{c}{2}),\\
&&\rho^{\pm}(u \pm \frac{c}{2})\Big( \bar{b}(u \pm \frac{c}{2},P+h)  L^{\pm}_{21}(u_1)  L^{\mp}_{11}(u_2)+\bar{c_0}(u \pm \frac{c}{2},P+h) L^{\pm}_{11}(u_1)  L^{\mp}_{21}(u_2)\Big)\nn\\
&&\qquad \quad  \qquad = L^{\mp}_{11}(u_2)  L^{\pm}_{21}(u_1) \rho^{\pm}(u \mp \frac{c}{2}),\nn\\
&& \rho^{\pm}(u \pm \frac{c}{2})\Big( \bar{b}(u \pm \frac{c}{2},P+h) L^{\pm}_{21}(u_1)  L^{\mp}_{12}(u_2)+\bar{c_0}(u \pm \frac{c}{2},P+h) L^{\pm}_{11}(u_1)  L^{\mp}_{22}(u_2)\Big)\nn\\
&&\qquad \quad   \qquad = \Big( L^{\mp}_{11}(u_2)  L^{\pm}_{22}(u_1) \bar{c_0}(u \mp \frac{c}{2},P) + L^{\mp}_{12}(u_2)  L^{\pm}_{21}(u_1) b(u \mp \frac{c}{2},P)  \Big)\ \rho^{\pm}(u \mp \frac{c}{2}),\nn\\ \\
&& \rho^{\pm}(u \pm \frac{c}{2})\Big(\bar{b}(u \pm \frac{c}{2},P+h)  L^{\pm}_{22}(u_1)  L^{\mp}_{11}(u_2)+\bar{c_0}(u \pm \frac{c}{2},P+h) L^{\pm}_{12}(u_1)  L^{\mp}_{21}(u_2) \Big) \nn\\
&&\qquad \quad  \qquad =\Big(L^{\mp}_{11}(u_2)  L^{\pm}_{22}(u_1) \bar{b}(u \mp \frac{c}{2},P) + L^{\mp}_{12}(u_2)  L^{\pm}_{21}(u_1) {c_0}(u \mp \frac{c}{2},P) \Big) \rho^{\pm}(u \mp \frac{c}{2}),\nn\\ 
&& \rho^{\pm}(u \pm \frac{c}{2}) \Big(\bar{b}(u \pm \frac{c}{2},P+h) L^{\pm}_{22}(u_1)  L^{\mp}_{12}(u_2)+\bar{c_0}(u \pm \frac{c}{2},P+h) L^{\pm}_{12}(u_1) L^{\mp}_{22}(u_2)\Big)\nn\\
&&\qquad \quad  \qquad = L^{\mp}_{12}(u_2)  L^{\pm}_{22}(u_1)\rho^{\pm}(u \mp \frac{c}{2}),
\ena
\bea
&&\rho^{\pm}(u \pm \frac{c}{2})  L^{\pm}_{21}(u_1)  L^{\mp}_{21}(u_2)= L^{\mp}_{21}(u_2)  L^{\pm}_{21}(u_1) \rho^{\pm}(u \mp \frac{c}{2}),\\
&&\rho^{\pm}(u \pm \frac{c}{2})  L^{\pm}_{21}(u_1)L^{\mp}_{22}(u_2)\nn\\
&&\qquad \quad  \qquad = \Big( L^{\mp}_{21}(u_2)  L^{\pm}_{21}(u_1) b(u \mp \frac{c}{2},P)   + \ L^{\mp}_{21}(u_2)  L^{\pm}_{22}(u_1) \bar{c_0}(u \mp \frac{c}{2},P) \Big) \rho^{\pm}(u \mp \frac{c}{2}),\nn\\ 
&&\rho^{\pm}(u \pm \frac{c}{2})  L^{\pm}_{22}(u_1)  L^{\mp}_{21}(u_2) \nn\\
&&\qquad \quad  \qquad = \Big(L^{\mp}_{21}(u_2)  L^{\pm}_{21}(u_1) {c_0}(u \mp \frac{c}{2},P) + \  L^{\mp}_{21}(u_2)  L^{\pm}_{22}(u_1) \bar{b}(u \mp \frac{c}{2},P)\Big) \ \rho^{\pm}(u \mp \frac{c}{2}),\nn\\ \lb{A.29} \\
&&\rho^{\pm}(u \pm \frac{c}{2})  L^{\pm}_{22}(u_1) L^{\mp}_{22}(u_2)= L^{\mp}_{22}(u_2)  L^{\pm}_{22}(u_1) \rho^{\pm}(u \mp \frac{c}{2}),
\ena
\end{prop}
where the dynamical $R$-matrix is expressed as
\begin{center}
\bea \lb{rappendix}
&&R^\pm(u,P)=\rho^\pm(u) \left(
\begin{array}{cccc}
1 & 0 & 0 & 0 \\
0 &b(u,P) & c(u,P) & 0 \\
0 & \bar{c}(u,P) & \bar{b}(u,P) & 0 \\
0 & 0 & 0 & 1 \\
\end{array}
\right).
\ena
\end{center}

\end{appendix}

\end{document}